\documentclass[draft,preprint,12pt,numbers,sort&compress]{elsarticle}
\usepackage{}
\usepackage{mathrsfs}
\usepackage{amssymb}
\usepackage{amsthm}
\usepackage{mathrsfs}
\usepackage[centertags]{amsmath}
\usepackage{amsfonts}

\usepackage{color}

\usepackage{ctex}
\usepackage{CJK}
\input amssym.def
\input amssym.tex

\date{}

\newtheorem{Theorem}{Theorem}[section]

\newtheorem{Lemma}{Lemma}[section]

\newcommand\R{\mbox{\bf R}}

\newcommand\SR{\mbox{\scriptsize\bf R}}

\newcommand{\definition}{{\lower .5ex
  \hbox{$\>\>\stackrel{\triangle}{=}\>\>$} }}
\newcommand\supp{\mathop{\rm supp}}

\begin{document}

\baselineskip=22pt
\thispagestyle{empty}

\mbox{}
\bigskip

\begin{center}
{\Large \bf  The Cauchy problem for the  generalized  Ostrovsky equation with negative dispersion}\\[1ex]

{Xiangqian Yan\footnote{Email:yanxiangqian213@126.com}$^a$,\quad  Wei Yan\footnote{ Email:011133@htu.edu.cn}$^a$,
\quad}\\[1ex]

{$^a$College of Mathematics and Information Science, Henan Normal University,}\\
{Xinxiang, Henan 453007,   China}\\[1ex]

\end{center}

\bigskip
\bigskip

\noindent{\bf Abstract.}  This paper is devoted to studying the Cauchy problem for the generalized
Ostrovsky equation
\begin{eqnarray*}
      u_{t}-\beta\partial_{x}^{3}u-\gamma\partial_{x}^{-1}u+\frac{1}{k+1}(u^{k+1})_{x}=0,k\geq5
    \end{eqnarray*}
    with $\beta\gamma<0,\gamma>0$.
Firstly,   we prove that the Cauchy problem for  the generalized Ostrovsky equation
 is locally well-posed in $H^{s}(\R)\left(s>\frac{1}{2}-\frac{2}{k}\right)$.
 Then, we prove that the Cauchy problem for  the generalized Ostrovsky equation
  is locally well-posed in
 $X_{s}(\R): =\|f\|_{H^{s}}+\left\|\mathscr{F}_{x}^{-1}\left(\frac{\mathscr{F}_{x}
 f(\xi)}{\xi}\right)\right\|_{H^{s}}\left(s>\frac{1}{2}-\frac{2}{k}\right).$
Finally,  we  show that  the solution to  the Cauchy problem for generalized  Ostrovsky  equation
converges to the solution to the generalized KdV equation
as the rotation parameter $\gamma$
tends to zero for data belonging to $X_{s}(\R)(s>\frac{3}{2})$.
 The main difficulty is that the phase function
of Ostrosvky equation with negative dispersive
 $\beta\xi^{3}+\frac{\gamma}{\xi}$ possesses the zero  singular point.

\noindent {\bf Keywords}: Generalized Ostrovsky equation; Multilinear estimates; Cauchy problem; Weak rotation limit

\bigskip
\noindent {\bf AMS  Subject Classification}:  35G25
\bigskip

\leftskip 0 true cm \rightskip 0 true cm

\newpage{}

\begin{center}

{\Large \bf The Cauchy problem for the generalized  Ostrovsky  equation with negative dispersion}\\[1ex]

{Xiangqian Yan\footnote{Email:yanxiangqian213@126.com}$^a$,\quad  Wei Yan\footnote{ Email:011133@htu.edu.cn}$^a$,
\quad}\\[1ex]

{$^a$College of Mathematics and Information Science, Henan Normal University,}\\
{Xinxiang, Henan 453007,   China}\\[1ex]

\end{center}
 This paper is devoted to studying the Cauchy problem for the generalized
Ostrovsky equation
\begin{eqnarray*}
      u_{t}-\beta\partial_{x}^{3}u-\gamma\partial_{x}^{-1}u+\frac{1}{k+1}(u^{k+1})_{x}=0,k\geq5
    \end{eqnarray*}
    with $\beta\gamma<0,\gamma>0$.
Firstly,   we prove that the Cauchy problem for  the generalized Ostrovsky equation
 is locally well-posed in $H^{s}(\R)\left(s>\frac{1}{2}-\frac{2}{k}\right)$.
 Then, we prove that the Cauchy problem for  the generalized Ostrovsky equation
  is locally well-posed in
 $X_{s}(\R): =\|f\|_{H^{s}}+\left\|\mathscr{F}_{x}^{-1}\left(\frac{\mathscr{F}_{x}
 f(\xi)}{\xi}\right)\right\|_{H^{s}}\left(s>\frac{1}{2}-\frac{2}{k}\right).$
Finally,  we  show that  the solution to  the Cauchy problem for generalized  Ostrovsky  equation
converges to the solution to the generalized KdV equation
as the rotation parameter $\gamma$
tends to zero for data belonging to $X_{s}(\R)(s>\frac{3}{2})$.
 The main difficulty is that the phase function
of Ostrosvky equation with negative dispersive
 $\beta\xi^{3}+\frac{\gamma}{\xi}$ possesses the zero  singular point.

\bigskip

{\large\bf 1. Introduction}
\bigskip

\setcounter{Theorem}{0} \setcounter{Lemma}{0}

\setcounter{section}{1}

\indent In this paper, we consider the Cauchy problem for the
generalized Ostrovsky equation
\begin{eqnarray}
     && u_{t}-\beta\partial_{x}^{3}u+\frac{1}{k+1}(u^{k+1})_{x}
     -\gamma \partial_{x}^{-1}u=0,\beta<0,\gamma>0,\label{1.01}\\
     &&u(x,0)=u_{0}(x).\label{1.02}
 \end{eqnarray}
Here $\partial_{x}^{-1}$  is defined by
\begin{eqnarray*}
&&\partial_{x}^{-1}f(x)=\mathscr{F}^{-1}((i\xi)^{-1}\mathscr{F}f(\xi))
=\frac{1}{2}\left(\int_{-\infty}^{x}f(y)dy-\int_{x}^{\infty}f(y)dy\right).
 \end{eqnarray*}

\indent This equation was introduced by Levandosky and Liu in \cite{LL}.
When $k=2$,  (\ref{1.01}) was Ostrovsky equation with negative dispersion,
 which was  introduced by Ostrovsky in \cite{O} as a
 model for weakly nonlinear long waves, by taking into account of the
  Coriolis force, to describe the propagation of surface waves in the
   ocean in a rotating frame of reference\cite{B,GS,GGS}. The Ostrovsky equation
    with negative dispersion
       has been investigated
      by some authors \cite{GH,HJ,IM2006,IM2007,IM2009,I,LM,LHY,T}.
 In the absence of rotation (that is,  $\gamma=0$),
       it becomes the generalized Korteweg-de Vries equation, which has been
        investigated by some authors
       \cite{KPV1991,KPVCPAM,F,FLP2011,FP2013,FLPV2018,FP2018,FHR2019}.
       Kenig et al. \cite{KPVCPAM} established the small
data global theory of generalized Korteweg-de Vries equation in
the critical Sobolev space generalized Korteweg-de Vries equation in
the critical Sobolev space $\dot{H}^{s_{k}}(\R)$ with $s_{k}=\frac{1}{2}-\frac{2}{k}.$
 Farah and Pastor
 \cite{FP2013} presented  an alternative proof of the Kenig, Ponce and Vega \cite{KPVCPAM}.

To the best of our knowledge,
the optimal regularity problem  of he Cauchy problem for  (\ref{1.01}) in
inhomogeneous Sobolev space has not been investigated.
In this paper, we investigate the Cauchy problem for (\ref{1.01}) in inhomogeneous Sobolev space.
By using the Fourier restriction norm method introduced in \cite{BGAFA-Sch,BGAFA-KdV} and
developed in \cite{KPV1993,KPV1996,KM}, firstly  we establish three  multilinear estimates.
 Then, by  exploiting the multilinear estimates and the fixed point theorem,
  we prove that the Cauchy problem for  the generalized Ostrovsky equation is
  locally well-posed in $H^{s}(\R)\left(s>\frac{1}{2}-\frac{2}{k}\right)$ and
  in $X_{s}(\R): =\|f\|_{H^{s}}+\left\|\mathscr{F}_{x}^{-1}\left(\frac{\mathscr{F}_{x}f(\xi)}
  {\xi}\right)\right\|_{H^{s}}\left(s>\frac{1}{2}-\frac{2}{k}\right).$
Finally,  we  show that  the solution to  the Cauchy problem for generalized  Ostrovsky  equation
converges to the solution to the generalized KdV equation
as the rotation parameter $\gamma$
tends to zero for data belonging to $X_{s}(\R)(s>\frac{3}{2})$.

 We give some notations before presenting the main results. $\chi_{A}(x)=1$ if $x\in A$, otherwise $\chi_{A}(x)=0$.
$a\sim b$ means that there exists two positive constants $C_{1},C_{2}$ which
may depend on $\beta,\gamma$
such that $C_{1}|a|\leq |b|\leq C_{2}|a|.$
$a\gg b$ means that there exists  positive constant $C$ which may depend on $\beta,\gamma$
such that $|a|\geq C |b|.$ We define $A:=  max\left\{1,\left|\frac{6\gamma}{7\beta}\right|^{\frac{1}{4}},
\left|\frac{\gamma}{3\beta}\right|^{\frac{1}{2}},\left|\frac{\gamma}{\beta}\right|,100|\beta|,
100|\gamma|\right\},a:=2^{[A]}$, where $[A]$ denotes the largest integer which
is smaller than $A$.
We define $\langle\cdot\rangle=1+|\cdot|.$ Let $\psi$ be a smooth jump
 function, $\psi_{\delta}(t)=\psi(\frac{t}{\delta})$, satisfying $0\leq\psi \leq1$,
 $\psi(t)=1$ for $|t|\leq 1$, $supp\psi\in[-2,2]$ and $\psi(t)=0$ for $|t|>2.$
\begin{align*}
&\phi(\xi)=\beta\xi^{3}+\frac{\gamma}{\xi},\sigma=\tau+\phi(\xi),\\
&\mathscr{F}_{x}f(\xi)=\frac{1}{\sqrt{2\pi}}\int_{\SR}e^{-ix \xi}f(x)dx.
\end{align*}
\begin{align*}
&\mathscr{F}_{x}^{-1}f(\xi)=\frac{1}{\sqrt{2\pi}}\int_{\SR}e^{ix \xi}f(x)dx,\\
&\mathscr{F}u(\xi,\tau)=\frac{1}{2\pi}\int_{\SR^{2}}e^{-ix\xi-it\tau}u(x,t)dxdt,\\
&\mathscr{F}^{-1}u(\xi,\tau)=\frac{1}{2\pi}\int_{\SR^{2}}e^{ix\xi+it\tau}u(x,t)dxdt,\\
&D_{x}^{\alpha}u_{0}=\frac{1}{\sqrt{2\pi}}\int_{\SR} |\xi|^{\alpha}e^{ix\xi}
\mathscr{F}_{x}{u}_{0}(\xi)d\xi,\\
&J_{x}^{\alpha}u_{0}=\frac{1}{\sqrt{2\pi}}\int_{\SR} \langle\xi\rangle^{\alpha}e^{ix\xi}
\mathscr{F}_{x}{u}_{0}(\xi)d\xi,\\
&U^{\gamma,\beta}(t)u_{0}=\frac{1}{\sqrt{2\pi}}\int_{\SR} e^{ix\xi-it\phi(\xi)}
\mathscr{F}_{x}{u}_{0}(\xi)d\xi,\\
&D_{t}^{\alpha}U^{\gamma,\beta}(t)u_{0}=\frac{1}{\sqrt{2\pi}}\int_{\SR}
\left|\phi(\xi)\right|^{\alpha}e^{ix\xi-it\phi(\xi)}
\mathscr{F}_{x}{u}_{0}(\xi)d\xi,\\
&J_{t}^{\alpha}U^{\gamma,\beta}(t)u_{0}=\frac{1}{\sqrt{2\pi}}\int_{\SR}
\left\langle\phi(\xi)\right\rangle^{\alpha}e^{ix\xi-it\phi(\xi)}
\mathscr{F}_{x}{u}_{0}(\xi)d\xi,\\
&\|f\|_{L_{t}^{p}L_{x}^{q}}=\left(\int_{\SR}
\left(\int_{\SR}|f(x,t)|^{q}dx\right)^{\frac{p}{q}}dt\right)^{\frac{1}{p}},\\
&P_{N}u=\frac{1}{\sqrt{2\pi}}\int_{|\xi|<N}e^{ix\xi}\mathscr{F}_{x}u(\xi)d\xi,\\
&P^{N}u=\frac{1}{\sqrt{2\pi}}\int_{|\xi|\geq N}e^{ix\xi}\mathscr{F}_{x}u(\xi)d\xi,\\
&\|f\|_{L_{xt}^{p}}=\|f\|_{L_{x}^{p}L_{t}^{p}}.
\end{align*}
$H^{s}(\R)=\left\{f\in \mathscr{S}^{\prime}(\R):\|f \|_{H^{s}(\SR)}=
 \|\langle\xi\rangle ^{s}\mathscr{F}_{x}{f}\|_{L_{\xi}^{2}(\SR)}<\infty\right\}$.
 The space $X_{s,b}(\R^{2})$ is defined as follows:
\begin{eqnarray*}
X_{s,b}(\R^{2})=\left\{u\in \mathscr{S}^{\prime}(\R^{2}):\|u\|_{X_{s,b}}
=\left[\int_{\SR^{2}}\langle \xi\rangle^{2s}\langle \sigma\rangle^{2b}
|\mathscr{F}u(\xi,\tau)|^{2}d\xi\right]^{\frac{1}{2}}<\infty\right\}.
\end{eqnarray*}
Here, $\langle \sigma\rangle =1+|\tau+\phi(\xi)|.$
The space $X_{s}(\R)$ is defined $X_{s}(\R)=\{f\in H^{s}(\R):\mathscr{F}_{x}^{-1}
(\frac{\mathscr{F}_{x}f(\xi)}{\xi})\in H^{s}(R)\}$ with the norm
 $\|f\|_{X_{s}(\SR)}=\|f\|_{H^{s}}+\left\|\mathscr{F}_{x}^{-1}
\left(\frac{\mathscr{F}_{x}f(\xi)}{\xi}\right)\right\|_{H^{s}}$.
The $\tilde{X}_{s,b}$ is defined as follows
\begin{eqnarray*}
&&\|u\|_{\tilde{X}_{s,b}}=\|u\|_{X_{s,b}}+\|\partial_{x}^{-1}u\|_{X_{s,b}}.
\end{eqnarray*}

The main result of this paper are as follows:
\begin{Theorem}\label{Thm1}
(\ref{1.01}) is locally well-posed for the initial data $u_{0}$ in $H^{s}(\R)$ with
 $s>\frac{1}{2}-\frac{2}{k},k\geq5$ and $\beta<0,\gamma>0$ .
\end{Theorem}

\noindent{\bf Remark1.} Lemma 3.1 in combination with Lemma 2.1 and the
 fixed point argument yields
Theorem 1.1. Thus, Lemma 3.1 plays the crucial role in establishing Theorem 1.1.
  The structure of (\ref{1.01})
 is more complicated than its of  the generalized KdV
equation. More precisely, the phase function of (\ref{1.01}) $\beta\xi^{3}
+\frac{\gamma}{\xi}$ possesses the zero singular point.
For the generalized KdV equation,  the following two important facts are valid.
\begin{eqnarray}
&\left\|\int_{\xi=\xi_{1}+\xi_{2},\tau=\tau_{1}+\tau_{2}}
\left|\xi_{1}^{2}-\xi_{2}^{2}\right|^{\frac{1}{2}}
\mathscr{F}{u_{1}}(\xi_{1},\tau_{1})\mathscr{F}{u_{2}}
(\xi_{2},\tau_{2})d\xi_{1}d\tau_{1}\right\|_{L_{\xi\tau}^{2}}\nonumber\\&\leq C\prod \limits_{j=1}^{2}\left\|(1+|\tau_{j}-\xi_{j}^{3}|)^{b}\mathscr{F}u_{j}\right\|_{L_{\xi_{j}\tau_{j}}^{2}}\label{1.03}
\end{eqnarray}
and
\begin{eqnarray}
\left\|P_{M}u\right\|_{L_{x}^{2}L_{t}^{\infty}}\leq C\prod \limits_{j=1}^{2}
\left\|(1+|\tau-\xi^{3}|)^{b}\mathscr{F}u\right\|_{L_{\xi\tau}^{2}}\label{1.04}
\end{eqnarray}
where $0<M\leq1.$
For the proof of (\ref{1.03}) and (\ref{1.04}),  we refer the readers to \cite{G,KPV1993}.

\noindent For (\ref{1.01}),  the following two important facts are valid.
\begin{eqnarray}
&\left\|\int_{\xi=\xi_{1}+\xi_{2},\tau=\tau_{1}+\tau_{2}}
\left|\phi^{\prime}(\xi_{1})-\phi^{\prime}(\xi_{2})\right|^{\frac{1}{2}}
\mathscr{F}{u_{1}}(\xi_{1},\tau_{1})\mathscr{F}{u_{2}}
(\xi_{2},\tau_{2})d\xi_{1}d\tau_{1}\right\|_{L_{\xi\tau}^{2}}\nonumber\\&\leq C\prod \limits_{j=1}^{2}\left\|(1+\left|\tau_{j}+\beta\xi_{j}^{3}+\frac{\gamma}{\xi_{j}}|\right)^{b}\mathscr{F}u_{j}
\right\|_{L_{\xi_{j}\tau_{j}}^{2}}\label{1.05}
\end{eqnarray}
and
\begin{eqnarray}
\left\|P_{M}u\right\|_{L_{x}^{2}L_{t}^{\infty}}\leq C\prod \limits_{j=1}^{2}\left\|D_{x}^{-s_{1}}(1+\left|\tau+\beta\xi_{j}^{3}+\frac{\gamma}{\xi_{j}}\right|)^{b}
\mathscr{F}u\right\|_{L_{\xi\tau}^{2}}\label{1.06}
\end{eqnarray}
for $0<M\leq1,s_{1}>\frac{1}{4}.$

\noindent For (\ref{1.05}),  only when
\begin{eqnarray}
(\xi_{1},\tau_{1},\xi,\tau)\in
\Omega :=\left\{(\xi_{1},\tau_{1},\xi,\tau)\in \R^{4}:
\left|1+\frac{\gamma}{3\beta\xi_{1}^{2}\xi_{2}^{2}}\right|\geq C\right\}\label{1.07}
\end{eqnarray}
for some $0<C<1$, then we have
\begin{eqnarray}
&\left\|\int_{\xi=\xi_{1}+\xi_{2},\tau=\tau_{1}+\tau_{2}}
\left|\xi_{1}^{2}-\xi_{2}^{2}\right|^{\frac{1}{2}}
\mathscr{F}{u_{1}}(\xi_{1},\tau_{1})\mathscr{F}{u_{2}}
(\xi_{2},\tau_{2})d\xi_{1}d\tau_{1}\right\|_{L_{\xi\tau}^{2}}
\nonumber\\&\leq C\prod \limits_{j=1}^{2}\left\|(1+
\left|\tau_{j}+\beta\xi_{j}^{3}+\frac{\gamma}{\xi_{j}}|\right)^{b}\mathscr{F}u_{j}
\right\|_{L_{\xi_{j}\tau_{j}}^{2}}\label{1.08}.
\end{eqnarray}
When
\begin{eqnarray}
(\xi_{1},\tau_{1},\xi,\tau)\in
\Omega^{c} :=\left\{(\xi_{1},\tau_{1},\xi,\tau)\in \R^{4}:
\left|1+\frac{\gamma}{3\beta\xi_{1}^{2}\xi_{2}^{2}}\right|<
 C\right\}\label{1.09},
\end{eqnarray}
 in most cases, we consider
\begin{eqnarray*}
\left|1+\frac{\gamma}{3\beta\xi^{2}\xi_{k+1}^{2}}\right|\geq\frac{1}{2},
\left|1+\frac{\gamma}{3\beta\xi^{2}\xi_{k+1}^{2}}\right|<\frac{1}{2},
\end{eqnarray*}
which are just (\ref{3.08}),  (\ref{3.09}).
When (\ref{3.09}) is valid, we will use (\ref{2.010}), (\ref{2.057}).
Moreover,  the maximal function estimate (\ref{2.045}) plays the important  role
in dealing with  $(\xi_{1},\xi_{2},\cdot\cdot\cdot,\xi_{k},\xi,\tau_{1},\tau_{2},
\cdot\cdot\cdot,\tau_{k},\tau)\in \Omega_{5}$.

\noindent{\bf Remark2.} From \cite{LM}, we know that when $u_{0}\in H^{1}(\R),$
 we can not obtain the upper bound of
$\|u\|_{H^{1}(\SR)}$, thus, we cannot obtain the global well-posedness of (\ref{1.01}).

\noindent{\bf Remark3.} From \cite{KPVCPAM} and
 \cite{FP2013}, we know that $s=\frac{1}{2}-\frac{2}{k}$ is the critical regularity
 indices in Sobolev spaces for (\ref{1.01}) with $\gamma=0$, which is just the generalized KdV equation.

\begin{Theorem}\label{Thm2}
(\ref{1.01}) is locally well-posed for the initial data $u_{0}$
 in $X_{s}(\R)$ with
 $s>\frac{1}{2}-\frac{2}{k},k\geq5$ and $\beta<0,\gamma>0$.
\end{Theorem}

\begin{Theorem}\label{Thm3}
Let $u^{\gamma}$ be the solution to (\ref{1.01}) in $X_{s}(\R)$
with $s>\frac{3}{2},k\geq5$
 and $\beta<0,\gamma>0$.Then $u^{\gamma}$ converges
  in $H^{s}(\R)$ with
  $s>\frac{3}{2},k\geq5$ to the solution to generalized
   KdV equation as $\gamma\rightarrow0$ and $u_{0}$ tends to the initial data of
    the generalized KdV in $L^{2}(\R)$.
   More precisely,
   \begin{eqnarray}
\|u-v\|_{L^{2}}\leq e^{CT\sup\limits_{t\in [0,T]}\left[\|u\|_{X_{s}}+\|v\|_{X_{s}}\right]^{k}}
\left[\|u_{0}-v_{0}\|_{L^{2}}+C|\gamma| T\sup\limits_{t\in [0,T]}
\|u\|_{X_{s}}\right].\label{1.010}
\end{eqnarray}
Here, $T$    is   the time  lifespan of the solution  to  (\ref{1.01})  for data
 in $X_{s}(\R)$ guaranteed by Theorem 1.2 and   $u$ is the solution to (\ref{1.01})-(\ref{1.02})
  and $v$ is the solution to the Cauchy problem for the  generalized KdV equation
\begin{eqnarray*}
     && v_{t}-\beta\partial_{x}^{3}v+\frac{1}{k+1}(v^{k+1})_{x}=0,\\
     &&v(x,0)=v_{0}(x).
 \end{eqnarray*}
\end{Theorem}

The rest of the paper is arranged as follows. In Section 2,
  we give some preliminaries.
In Section 3, we show some multilinear estimates. In Section 4,
 we prove Theorem 1.1.
In Section 5, we prove Theorem 1.2. In Section 6, we prove Theorem 1.3.

\bigskip

 \noindent{\large\bf 2. Preliminaries }

\setcounter{equation}{0}

\setcounter{Theorem}{0}

\setcounter{Lemma}{0}

\setcounter{section}{2}

\begin{Lemma}\label{Lemma2.1}
Let $\delta\in (0,1)$ and $s\in \R$ and $-\frac{1}{2}<b^{\prime}
\leq0\leq b\leq b^{\prime}+1$.
Then, for $h\in X_{s,b^{\prime}},$  we have
\begin{eqnarray}
&&\left\|\psi(t)U^{\gamma,\beta}(t)h\right\|_{X_{s,\frac{1}{2}+\epsilon}}\leq C\|h\|_{H^{s}},\label{2.01}\\
&&\left\|\psi\left(\frac{t}{\delta}\right)\int_{0}^{t}U^{\gamma,\beta}(t-\tau)h(\tau)
d\tau\right\|_{X_{s,b}}\leq C
\delta^{1+b^{\prime}-b}\|h\|_{X_{s,b^{\prime}}}.\label{2.02}
\end{eqnarray}
\end{Lemma}

For the proof of Lemma 2.1, we refer the readers to \cite{BGAFA-Sch,BGAFA-KdV,G}.

\begin{Lemma}\label{Lemma2.2}
Let $q\geq8, 0<M\leq 1$ and $s=\frac{1}{8}-\frac{1}{q}$,
$s_{1}=\frac{1}{4}+\epsilon$
and $b=\frac{1}{2}+\frac{\epsilon}{24}$
 and $0\leq\epsilon\leq 10^{-3}$. Then,  we have
\begin{align}
&\|U^{\gamma,\beta}(t)u_{0}\|_{L_{xt}^{8}}\leq C\|u_{0}\|_{L_{x}^{2}},\label{2.03}\\
&\|u\|_{L_{xt}^{q}}\leq C\|u\|_{X_{4s,b}},\label{2.04}\\
&\|D_{x}^{\frac{1}{6}}P^{a}u\|_{L_{xt}^{6}}\leq C\|u\|_{X_{0,b}},\label{2.05}\\
&\|u\|_{L_{xt}^{8}}\leq C\|u\|_{X_{0,b}}\label{2.06},\\
&\|u\|_{L_{xt}^{\frac{8}{1+\epsilon}}}\leq C\|u\|_{X_{0,\frac{3-\epsilon}{3}b}}
\leq C\|u\|_{X_{0,\frac{1}{2}-\frac{\epsilon}{12}}}\label{2.07},\\
&\|D_{x}P^{a}u\|_{L_{x}^{\infty}L_{t}^{2}}\leq C\|u\|_{X_{0,b}},\label{2.08}\\
&\|D_{x}^{s_{1}}\psi \left(t\right)P_{M}u\|_{L_{x}^{2}L_{t}^{\infty}}\leq C\|u\|_{X_{0,b}},\label{2.09}\\
&\|D_{x}^{1-2\epsilon}P^{a}u\|_{L_{x}^{\frac{1}{\epsilon}}L_{t}^{2}}\leq C
\|u\|_{X_{0,(1-2\epsilon)b}}\label{2.010},\\
&\|D_{x}^{\frac{3-\epsilon}{18}}P^{a}u\|_{L_{xt}^{\frac{6}{1+\epsilon}}}
\leq C\|u\|_{X_{0,\frac{3-\epsilon}{3}b}}\label{2.0888888}.
\end{align}
\end{Lemma}
\noindent {\bf Proof.}For (\ref{2.03}),   we  refer  the readers to  Lemma 2.1 of \cite{WC}.
 By using the Sobolev embeddings theorem and (\ref{2.03}), we derive
\begin{eqnarray}
&&\left\|U^{\gamma,\beta}(t)u_{0}\right\|_{L_{xt}^{q}}\leq C\left\|D_{x}^{s}D_{t}^{s}U^{\gamma,\beta}(t)u_{0}\right\|_{L_{xt}^{8}}
= C\left\|\int_{\SR}e^{ix\xi-it\phi}\left|\beta\xi^{3}+\frac{\gamma}{\xi}\right|^{s}
|\xi|^{s}\mathscr{F}_{x}u_{0}(\xi)d\xi\right\|_{L_{xt}^{8}}
\nonumber\\
&&\leq C\left\||\xi^4+\frac{\gamma}{\beta}|^{s}\mathscr{F}_{x}u_{0}\right\|_{L_{\xi}^{2}}\leq C\left\||\xi^{2}+a|^{2s}\mathscr{F}_{x}u_{0}\right\|_{L_{\xi}^{2}}\nonumber\\
&&\leq C\left\|\xi^{4s}u_{0}\right\|_{L_{\xi}^{2}}+C
\left\|u_{0}\right\|_{L_{\xi}^{2}}\leq C\|u_{0}\|_{H^{4s}}
\label{2.011}.
\end{eqnarray}
By changing variable $\tau=\lambda-\phi(\xi)$, we derive
 \begin{eqnarray}
&&u(x,t)=\frac{1}{2\pi}\int_{\SR^{2}}e^{ix\xi+it\tau}\mathscr{F}u(\xi,\tau)d\xi d\tau\nonumber\\
&&=\frac{1}{2\pi}\int_{\SR^{2}}e^{ix\xi+it(\lambda-\phi(\xi))}\mathscr{F}u(\xi,\lambda-\phi(\xi))d\xi d\lambda\nonumber\\
&&=\frac{1}{2\pi}\int_{\SR}e^{it\lambda}\left(\int_{\SR}e^{ix\xi-it\phi(\xi)}
\mathscr{F}u(\xi,\lambda-\phi(\xi)) d\xi\right)d\lambda\label{2.012}.
 \end{eqnarray}
 By using  (\ref{2.011}),  (\ref{2.012}) and Minkowski's inequality,  for  $b>\frac{1}{2},$  we derive
\begin{eqnarray}
&&\left\|u\right\|_{L_{xt}^{q}}\leq C\int_{\SR}\left\|\left(\int_{\SR}e^{ix\xi-it\phi(\xi)}
\mathscr{F}u(\xi,\lambda-\phi(\xi))d\xi\right)\right\|_{L_{xt}^{q}} d\lambda\nonumber\\&&\leq C\int_{\SR}\left\|\mathscr{F}u(\xi,\lambda-\phi(\xi))\right\|_{H^{4s}}d\lambda\nonumber\\&&\leq C
\left[\int_{\SR}(1+|\lambda|)^{2b}\left\|\mathscr{F}u(\xi,\lambda-\phi(\xi))\right\|_{H^{4s}}^{2}d\lambda\right]^{\frac{1}{2}}
\left[\int_{\SR}(1+|\lambda|)^{-2b}d\lambda\right]^{\frac{1}{2}}\nonumber\\
&&\leq C\left[\int_{\SR}(1+|\tau+\phi(\xi)|)^{2b}\left\|\mathscr{F}u(\xi,\tau)\right\|_{H^{4s}}^{2}d\tau\right]^{\frac{1}{2}}
=\left\|u\right\|_{X_{4s,b}}.\label{2.013}
\end{eqnarray}
Thus, we obtain (\ref{2.04}).
 \begin{eqnarray}
&\|D_{x}^{\frac{1}{6}}U^{\gamma,\beta}(t)P^{a}u_{0}\|_{L_{xt}^{6}}\leq C\|u_{0}\|_{L^{2}}\label{2.014}
\end{eqnarray}
can be seen in (2.3) of Lemma 2.1 of \cite{GH}.
 By using (\ref{2.014}) and a proof similar to (\ref{2.04}), we obtain (\ref{2.05}).
 By using (\ref{2.03}) and a proof similar to (\ref{2.04}), we obtain (\ref{2.06}).
 Interpolating (\ref{2.06}) with
\begin{eqnarray}
&&\left\|u\right\|_{L_{xt}^{2}}\leq C\left\|u\right\|_{X_{0,0}}\label{2.015}
\end{eqnarray}
yields (\ref{2.07}).
\begin{eqnarray}
&\|D_{x}U^{\gamma,\beta}(t)P^{a}u_{0}\|_{L_{x}^{\infty}L_{t}^{2}}\leq C\|u_{0}\|_{L_{x}^{2}}\label{2.016}
\end{eqnarray}
can be seen in (2.1) of Lemma 2.1 of \cite{GH}.
By using (\ref{2.016}) and a proof similar to (\ref{2.04}), we obtain (\ref{2.08}).
From (31) of Lemma 3.4 of \cite{LM}, for $s_{1}=\frac{1}{4}+\epsilon,
s_{2}>\frac{3}{4},0<M\leq 1$,  we know that
\begin{eqnarray}
&\|U^{\gamma,\beta}(t)P_{M}\psi(t)u_{0}\|_{L_{x}^{2}L_{t}^{\infty}}\leq C
\|D_{x}^{-s_{1}}P_{M}u_{0}\|_{L^{2}}+C\|P_{M}u_{0}\|_{H^{s_{2}}}
\nonumber\\&\leq C\|D_{x}^{-s_{1}}P_{M}u_{0}\|_{L^{2}}\label{2.017}.
\end{eqnarray}
By using (\ref{2.017}) and a proof similar to (\ref{2.04}), we obtain
(\ref{2.09}). Interpolating (\ref{2.08}) with (\ref{2.015}) leads to  (\ref{2.010}).
Interpolating (\ref{2.05}) with (\ref{2.015}) leads to  (\ref{2.0888888}).

We have completed the proof of Lemma 2.2.

\begin{Lemma}\label{Lemma2.3}
We assume that $f\in C_{0}(\R)$ and $\phi\in C^{\infty}(\R)$
 and $x_{i}(1\leq i\leq n)$
 are only simple zeros of  $\phi(x)$ which means $\phi(x_{i})=0$
  and $\phi^{\prime}(x_{i})\neq0$. Then we have
\begin{eqnarray}
&&\delta[\phi(x)]=\sum\limits_{i=1}^{n}\frac{\delta(x-x_{i})}{|\phi^{\prime}(x_{i})|}
,\label{2.018}
\end{eqnarray}
where $\delta$ is Dirac delta function.
\end{Lemma}
\noindent {\bf Proof.} Since  $f\in C_{0}(\R)$ and $\phi\in C^{\infty}(\R)$
and $\phi^{\prime}(x_{i})\neq0$,
$\forall \epsilon >0(<\frac{|\phi^{\prime}(x_{i})|}{2}),$ there exists
$\delta_{1} >0$ such that when $|x-x_{i}|<\delta_{1}$, we have
\begin{eqnarray}
&&|f(x)-f(x_{i})|<\epsilon,\label{2.019}\\
&&|\phi^{\prime}(x)-\phi^{\prime}(x_{i})|<\epsilon.\label{2.020}
\end{eqnarray}
From (\ref{2.020}),  we know that
$|\phi^{\prime}(x)|\geq|\phi^{\prime}(x_{i})|-\epsilon\geq \frac{|\phi^{\prime}(x_{i})|}{2}.$
Thus, when $x\in (x_{i}-\delta_{1},x_{i}+\delta_{1})$, we
have
\begin{eqnarray}
&&\left|\frac{f(x)}{\phi^{\prime}(x)}-\frac{f(x_{i})}{\phi^{\prime}(x_{i})}\right|=
\left|\frac{f(x)\phi^{\prime}(x_{i})-f(x_{i})\phi^{\prime}(x)}{\phi^{\prime}(x_{i})\phi^{\prime}(x)}\right|\nonumber\\
&&\leq \left|\frac{|f(x)-f(x_{i})}{\phi^{\prime}(x)}\right|+\left|\frac{f(x_{i})
(\phi^{\prime}(x)-\phi^{\prime}(x_{i}))}{\phi^{\prime}(x)\phi^{\prime}(x_{i})}\right|\nonumber\\
&&\leq \frac{2\epsilon}{|\phi^{\prime}(x_{i})|}+\frac{2\epsilon|f(x_{i})|}{|\phi^{\prime}(x_{i})|^{2}}.\label{2.021}
\end{eqnarray}
We defined
\begin{eqnarray}
&&F_{i}=\int_{x_{i}-\delta_{1}}^{x_{i}+\delta_{1}}f(x)\delta[\phi(x)]dx.\label{2.022}
\end{eqnarray}
We claim
\begin{eqnarray}
F_{i}=\frac{f(x_{i})}{|\phi^{\prime}(x_{i})|}.\label{2.023}
\end{eqnarray}
By a change of the variable $u=\phi(x)$, $du=\phi^{\prime}(x)dx$. Then, by using (\ref{2.023})  and  $\int_{\phi(x_{i}-\delta_{1})}^{\phi(x_{i}+\delta_{1})}\delta(u)du=1$,  we have
\begin{eqnarray}
&&\left|F_{i}-\frac{f(x_{i})}{|\phi^{\prime}(x_{i})|}\right|=
\left|\int_{\phi(x_{i}-\delta_{1})}^{\phi(x_{i}+\delta_{1})}f(x)\delta(u)\frac{du}
{|\phi^{\prime}(x)|}-\int_{\phi(x_{i}-\delta_{1})}^{\phi(x_{i}+\delta_{1})}f(x_{i})\delta(u)\frac{du}
{|\phi^{\prime}(x_{i})|}\right|\nonumber\\
&&=\left|\int_{\phi(x_{i}-\delta_{1})}^{\phi(x_{i}+\delta_{1})}
\left|\frac{f(x)}{\phi^{\prime}(x)}-\frac{f(x_{i})}{\phi^{\prime}(x_{i})}\right|\delta(u)du\right|\nonumber\\
&&\leq \left[\frac{2\epsilon}{|\phi^{\prime}(x_{i})|}+\frac{2\epsilon|f(x_{i})|}{|\phi^{\prime}(x_{i})|^{2}}\right]
\int_{\phi(x_{i}-\delta_{1})}^{\phi(x_{i}+\delta_{1})}
\delta(u)du=\left[\frac{2\epsilon}{|\phi^{\prime}(x_{i})|}+
\frac{2\epsilon|f(x_{i})|}{|\phi^{\prime}(x_{i})|^{2}}\right].\label{2.024}
\end{eqnarray}
Since $\epsilon>0$ is arbitrary small, thus, from (\ref{2.024}),
we know  that (\ref{2.023}) is valid.
Using (\ref{2.023}),  we have
\begin{eqnarray}
&&\int_{-\infty}^{+\infty}f(x)\delta[\phi(x)]dx=\sum\limits_{i=1}^{n}F_{i}
=\sum\limits_{i=1}^{n}\frac{f(x_{i})}{|\phi^{\prime}(x_{i})|}\nonumber\\
&&=\int_{-\infty}^{+\infty}f(x)\sum\limits_{i=1}^{n}\left[\frac{\delta(x-x_{i})}
{|\phi^{\prime}(x_{i})|}\right]dx.\label{2.025}
\end{eqnarray}
Thus, since $f\in C_{0}(\R),$    we obtain
\begin{eqnarray}
&&\delta[\phi(x)]=\sum\limits_{i=1}^{n}\left[\frac{\delta(x-x_{i})}
{|\phi^{\prime}(x_{i})|}\right].\label{2.026}
\end{eqnarray}

We have completed the proof of Lemma 2.3.

\noindent{\bf Remark3.}  In page 184 of \cite{GS}, Gel'fand and Shilov \cite{GS} presented the
definition of (\ref{2.018}), however, they did not give the strict proof.

\begin{Lemma}\label{Lemma2.4}
\noindent Let $b>\frac{1}{2}$.
Then we have
\begin{eqnarray}
\left\|I^{\frac{1}{2}}(e^{t(\beta \partial_{x}^{3}+\partial_{x}^{-1}\gamma)}f_{1},
e^{t(\beta \partial_{x}^{3}+\partial_{x}^{-1}\gamma)}f_{2})\right\|
  _{L_{xt}^{2}}\leq C\|f_{1}\|_{L_{x}^{2}}\|f_{2}\|_{L_{x}^{2}}.\label{2.027}
\end{eqnarray}
Here
\begin{eqnarray}
\mathscr{F}I^{s}(f_{1},f_{2})(\xi,\tau)
&=&\int_{\!\!\!\mbox{\scriptsize $
\begin{array}{l}
\xi=\xi_{1}+\xi_{2}\\
\end{array}
$}}
\left|\phi^{\prime}(\xi_{1})-\phi^{\prime}(\xi_{2})\right|^{s}
\mathscr{F}_{x}{f_{1}}(\xi_{1})\mathscr{F}_{x}{f_{2}}
(\xi_{2})\,d\xi_{1}.\label{2.028}
\end{eqnarray}
In particular, for $b>\frac{1}{2},$   we have
\begin{eqnarray}
\left\|I^{\frac{1}{2}}(u_{1},u_{2})\right\|
  _{L_{xt}^{2}}\leq C\prod_{j=1}^{2}\| u_{j}\| _{X_{0,b}}.\label{2.029}
\end{eqnarray}
\end{Lemma}
\noindent {\bf Proof.}Followed   the idea of \cite{G,G2004}, we present the proof of Lemma 2.4.
By  using the  Plancherel identity with respect to  the space variable and
\begin{eqnarray*}
&&\int_{\SR}e^{-it(\phi(\xi_{1})+\phi(\xi_{2})-\phi(\eta_{1})-\phi(\eta_{2}))}dt\nonumber\\&&
=C\delta(\phi(\xi_{1})+\phi(\xi_{2})-\phi(\eta_{1})-\phi(\eta_{2})),
\end{eqnarray*}
we derive
\begin{eqnarray}
&&\left\|I^{\frac{1}{2}}\left(e^{t(\beta \partial_{x}^{3}+\partial_{x}^{-1}\gamma)}f_{1},e^{t(\beta \partial_{x}^{3}+\partial_{x}^{-1}\gamma)}f_{2}\right)\right\|_{L_{xt}^{2}}^{2}\nonumber\\&&
=C\int_{\SR^{2}}\left|\int_{\xi=\xi_{1}+\xi_{2}}
\left|\phi^{\prime}(\xi_{1})-\phi^{\prime}(\xi_{2})\right|^{\frac{1}{2}}
e^{-it(\phi(\xi_{1})+\phi(\xi_{2}))}\mathscr{F}_{x}{f_{1}}(\xi_{1})
\mathscr{F}_{x}{f_{2}}(\xi_{2})d\xi_{1}\right|^{2}dtd\xi\nonumber\\
&&=C\int_{\SR^{2}}\int_{\xi=\xi_{1}+\xi_{2}}\int_{\xi=\eta_{1}+\eta_{2}}e^{-it(\phi(\xi_{1})+\phi(\xi_{2})-\phi(\eta_{1})-\phi(\eta_{2}))}
\left(\left|\phi^{\prime}(\xi_{1})-\phi^{\prime}(\xi_{2})\right|
\left|\phi^{\prime}(\eta_{1})-\phi^{\prime}(\eta_{2})\right|\right)^{\frac{1}{2}}
\nonumber\\&&\qquad\qquad\times\prod_{i=1}^{2}\mathscr{F}_{x}{f_{i}}(\xi)
\overline{\mathscr{F}_{x}{f_{i}}(\eta_{i})}d\xi_{1}d\eta_{1}d\xi dt\nonumber\\
&&=C\int_{\SR}
\int_{\xi=\xi_{1}+\xi_{2}}\int_{\xi=\eta_{1}+\eta_{2}}\delta(\phi(\xi_{1})
+\phi(\xi_{2})-\phi(\eta_{1})-\phi(\eta_{2}))
\left(\left|\phi^{\prime}(\xi_{1})-\phi^{\prime}(\xi_{2})\right|
\left|\phi^{\prime}(\eta_{1})-\phi^{\prime}(\eta_{2})\right|\right)^{\frac{1}{2}}
\nonumber\\&&\qquad\qquad\times\prod_{i=1}^{2}\mathscr{F}_{x}{f_{i}}(\xi)
\overline{\mathscr{F}_{x}{f_{i}}(\eta_{i})}d\xi_{1}d\eta_{1}d\xi.\label{2.030}
\end{eqnarray}
 From Lemma  2.2,  we have $\delta[g(x)]=\sum\limits_{i=1}^{n}
 \frac{\delta(x-x_{i})}{|g^{\prime}(x_{i})|}$, where
\begin{eqnarray}
&&g(x)=3\xi\beta(x^{2}-\xi_{1}^{2}+\xi(\xi_{1}-x))-\frac{\gamma\xi}{x(\xi-x)}
+\frac{\gamma\xi}{\xi_{1}(\xi-\xi_{1})},\label{2.031}
\end{eqnarray}
since $g^{\prime\prime}\neq0$, then $g$ has only two simple zeros, $x_{1}=\xi_{1}$
 and $x_{2}=\xi-\xi_{1}$.\\
Hence $g^{\prime}(x_{1})=\phi^{\prime}(\xi_{2})-\phi^{\prime}(\xi_{1})$
 and $g^{\prime}(x_{2})=\phi^{\prime}(\xi_{1})-\phi^{\prime}(\xi_{2})$. Thus (\ref{2.030})
   can be rewritten as
\begin{eqnarray}
&&\hspace{-0.8cm}C\int_{\SR} \int_{\xi=\xi_{1}+\xi_{2}}\int_{\xi=\eta_{1}+\eta_{2}}\frac{1}{|\phi^{\prime}(\xi_{2})-\phi^{\prime}(\xi_{1})|}
\delta(\eta_{1}-\xi_{1})
\left(\left|\phi^{\prime}(\xi_{1})-\phi^{\prime}(\xi_{2})\right|
\left|\phi^{\prime}(\eta_{1})-\phi^{\prime}(\eta_{2})\right|\right)^{\frac{1}{2}}
\nonumber\\&&\qquad\qquad\times\prod_{i=1}^{2}\mathscr{F}_{x}{f_{i}}(\xi)
\overline{\mathscr{F}_{x}{f_{i}}(\eta_{i})}d\xi_{1}d\eta_{1}d\xi\nonumber\\
&&\hspace{-0.8cm}+C\int_{\SR} \int_{\xi=\xi_{1}+\xi_{2}}\int_{\xi=\eta_{1}+\eta_{2}}\frac{1}
{|\phi^{\prime}(\xi_{2})-\phi^{\prime}(\xi_{1})|}\delta(\eta_{1}-(\xi-\xi_{1}))
\left(\left|\phi^{\prime}(\xi_{1})-\phi^{\prime}(\xi_{2})\right|\left|\phi^{\prime}
(\eta_{1})-\phi^{\prime}(\eta_{2})\right|\right)^{\frac{1}{2}}
\nonumber\\&&\qquad\qquad\times\prod_{i=1}^{2}\mathscr{F}_{x}{f_{i}}(\xi)
\overline{\mathscr{F}_{x}{f_{i}}(\eta_{i})}d\xi_{1}d\eta_{1}d\xi\nonumber\\
&&=C\int_{\SR}\int_{\xi=\xi_{1}+\xi_{2}}\prod_{i=1}^{2}
|\mathscr{F}_{x}{f_{i}}(\xi_{i})|^{2}d\xi_{1}d\xi\nonumber\\&&\qquad\qquad
+C\int_{\SR}\int_{\xi=\xi_{1}+\xi_{2}}\mathscr{F}_{x}{f_{1}}(\xi_{1})
\overline{\mathscr{F}_{x}{f_{1}}(\xi_{2})}\mathscr{F}_{x}{f_{2}}(\xi_{2})
\overline{\mathscr{F}_{x}{f_{2}}(\xi_{1})}d\xi_{1}d\xi\nonumber\\
&&\leq C\left(\prod_{i=1}^{2}\|f_{i}\|_{L_{x}^{2}}^{2}+\|\mathscr{F}_{x}{f_{1}}
\mathscr{F}_{x}{f_{2}}\|_{L_{\xi}^{1}}^{2}\right)\leq C\prod_{i=1}^{2}\|f_{i}\|_{L_{x}^{2}}^{2}.
\label{2.032}
\end{eqnarray}
We define $\mathscr{F}v_{j\tau}(\xi)=\mathscr{F}u_{j}(\xi,\tau-\phi(\xi))(j=1,2).$
Thus, by using (\ref{2.032}), we have
\begin{eqnarray}
I^{\frac{1}{2}}(u_{1},u_{2})(x,t)
&=&\int_{\SR^{2}}I^{\frac{1}{2}}(e^{it\tau_{1}}v_{1\tau_{1}},e^{it\tau_{1}}v_{2\tau_{2}})(x,t)
     d\tau_{1}d\tau_{2}.\label{2.033}
\end{eqnarray}
By using a  direct computation, we have
\begin{eqnarray}
\mathscr{F}\left[I^{\frac{1}{2}}(e^{it\tau_{1}}v_{1\tau_{1}},e^{it\tau_{2}}v_{2\tau_{2}})\right](\xi,\tau)=
\mathscr{F}\left[I^{\frac{1}{2}}(u_{1},u_{2})\right](\xi,\lambda-\tau_{1}-\tau_{2}).\label{2.034}
\end{eqnarray}
By using the Minkowski's  inequality,  the Plancherel identity,  (\ref{2.033})  and  (\ref{2.034}), we have
\begin{eqnarray}
&&\left\|I^{\frac{1}{2}}(u_{1},u_{2})\right\|
  _{L_{xt}^{2}}\leq C\int_{\SR^{2}}\left\|I^{\frac{1}{2}}(e^{it\tau_{1}}v_{1\tau_{1}},e^{it\tau_{2}}v_{2\tau_{2}})
  \right\|_{L_{xt}^{2}}d\tau_{1}d\tau_{2}\nonumber\\
  &&\leq C\int_{\SR^{2}}\left\|I^{\frac{1}{2}}(v_{1\tau_{1}},v_{2\tau_{2}})
  \right\|_{L_{xt}^{2}}d\tau_{1}d\tau_{2}\nonumber\\
  &&\leq C\int_{\SR^{2}}\|v_{1\tau_{1}}\|_{L_{x}^{2}}\|v_{2\tau_{1}}
  \|_{L_{x}^{2}}d\tau_{1}d\tau_{2}\nonumber\\
  &&\leq C\int_{\SR}\|\mathscr{F}_{x}v_{1\tau_{1}}\|_{L_{\xi}^{2}}d\tau_{1}
  \int_{\SR}\|\mathscr{F}_{x}v_{2\tau_{2}}\|_{L_{\xi}^{2}}d\tau_{2}\nonumber\\
  &&\leq C\prod\limits_{j=1}^{2}\left\|\langle \tau_{j}\rangle ^{b}
  \mathscr{F}u_{j}(\xi_{j},\tau_{j}-\phi(\xi_{j})\right\|_{L_{\xi\tau}^{2}}
  \leq C\prod\limits_{j=1}^{2}\left\|u_{j}\right\|_{X_{0,b}}.\label{2.035}
\end{eqnarray}

We have completed the proof of Lemma 2.4.
\begin{Lemma}\label{Lemma2.5}
We assume that $b>\frac{1}{2}$, $0\leq s\leq \frac{1}{2}.$
Then, we have
\begin{eqnarray}
\left\|I^{s}(u_{1},u_{2})\right\|
  _{L_{xt}^{2}}\leq C\prod\limits_{j=1}^{2}\|u_{j}\| _{X_{0,\frac{2+2s}{3}b}}.\label{2.036}
\end{eqnarray}
\end{Lemma}

Lemma 2.5 can be proved similarly to Corollary 3.2 of \cite{G2004} and Lemma 2.2
 of \cite{LLY} with the aid of Lemma 2.4.

\begin{Lemma}\label{lem2.6}
Let  $\phi_{j}(j=1,2)\in C^{\infty}(\R)$  and $\supp \phi_{2}\subset (a,b)$.
If $\phi_{1}^{\prime}(\xi)\neq 0$ for all $\xi\in [a,b]$, then
\begin{align*}
\left|\int_{a}^{b}e^{i\lambda\phi_{1}(\xi)}\phi_{2}(\xi)d\xi\right|
\leq \frac{C}{|\lambda|^{k}}
\end{align*}
for all $k\geq 0$,  where the constant $C$ depends on $\phi_{1},\phi_{2}$
 and k.
\end{Lemma}

Lemma 2.6 can be seen in \cite{Stein-new}.
\begin{Lemma}\label{lem2.7}
Let $\phi_{4}\in C^{\infty}(\R)$ and  $\phi_{3}\in C^{3}(\R)$ and   $\supp\phi_{3}\subset (a,b)$ and
 $\left|\phi_{3}^{(3)}\right|\geq 1$ uniformly with respect to $\xi$. Then, we have
\begin{align*}
\left|\int_{a}^{b}e^{i\lambda\phi_{3}(\xi)}\phi_{4}(\xi)d\xi\right|\leq
 \frac{C}{|\lambda|^{\frac{1}{3}}}\left(|\phi_{4}(b)|+\int_{a}^{b}
 |\phi_{4}^{'}(\xi)|d\xi\right).
\end{align*}
\end{Lemma}

Lemma 2.7 can be seen in \cite{Stein-new}.

\begin{Lemma}\label{lem2.8}Let
\begin{align*}
 K(x,t)=\int_{\SR}e^{-it\phi(\xi)+ix\xi}\chi_{[N,4N]}(|\xi|)d\xi.
 \end{align*}
 Here, $N\in 2^{Z},N\geq a$.
For $\gamma\geq7$, we obtain
\begin{align}\label{2.037}
 \|K\|_{L_{x}^{\frac{\gamma}{2}}L_{t}^{\infty}}\leq C N^{\frac{\gamma-2}{\gamma}}.
 \end{align}
\end{Lemma}
\noindent{\bf Proof.} We define
\begin{eqnarray*}
&&\Omega_{1}:=\left\{(x,t)\in \R\times \R: |x|\leq \frac{1}{N}\right\},\\
&&\Omega_{2}:= \left\{(x,t)\in \R\times \R: |x|> \frac{1}{N},|x|\geq 4000aN^{2}t\right\},\\
&&\Omega_{3}:= \left\{(x,t)\in \R\times \R: |x|>\frac{1}{N},|x|< 4000aN^{2}t\right\}.
\end{eqnarray*}
Obviously,  $\R_{x}\times \R_{t}=\bigcup\limits_{j=1}^{3}\Omega_{j}$.
We define $\Omega_{x,i}:=\left\{t\in \R|(x,t)\in \Omega_{i}\right\}$ for a fixed $x\in \R$.
 Without loss of generality, we assume that $N\leq \xi\leq 4N.$
Assume that  $\eta=\frac{\xi}{N}$, then we have
\begin{align*}
 K(x,t)=N\int_{1}^{4}e^{-iN^{3}\eta^{3}t-\frac{i\gamma t}{N\eta}+ixN\eta}d\eta.
 \end{align*}
 By using a  direct computation and from the definition $\Omega_{x1}$, we have
 \begin{align}\label{2.038}
\left[\int_{|x|\leq \frac{1}{N}}\left[\sup\limits_{t\in \Omega_{x,1}}|K(x,t)|
\right]^{\frac{\gamma}{2}}dx\right]^{\frac{2}{\gamma}}\leq C
N\left[\int_{|x|\leq \frac{1}{N}} dx\right]^{\frac{2}{\gamma}}
\leq C N^{\frac{\gamma-2}{\gamma}}.
\end{align}
By using a  direct computation, we have $-it\phi(\xi)+ix\xi=-i\beta N^{3}\eta^{3}t-\frac{i\gamma t}{N\eta}+ixN\eta
=ixN(\eta-\frac{\beta N^{2}\eta^{3}t}{x}-\frac{\gamma t}{N^{2}x\eta}):=
ixN\phi_{5}(\eta),$ where $\phi_{5}(\eta)=\eta-\frac{\beta N^{2}\eta^{3}t}{x}
-\frac{\gamma t}{N^{2}x\eta}$. Obviously,
\begin{align*}
 |\phi_{5}^{\prime}(\eta)|=\left|1-\frac{3\beta N^{2}t\eta^{2}}{x}-
 \frac{\gamma t}{N^{2}x\eta^{2}}\right|\geq 1-48\left|
 \frac{t\beta N^{2}}{x}\right|-\left|\frac{\gamma t}{N^{2}x}\right|\geq 1-\frac{1}{8}-
 \frac{1}{8}\geq\frac{1}{2}
 \end{align*}
 for any $(x,t)\in \Omega_{2}$. Therefore, $\phi_{5}^{\prime}\neq0$ in this region.
  From Lemma 2.6,  we know that  $|K(x,t)|\leq C N(N|x|)^{-2}=N^{-1}x^{-2}$. Thus, we derive
 \begin{align}\label{2.039}
\left[\int_{|x|> \frac{1}{N}}\left[\sup\limits_{t\in \Omega_{x,2}}|K(x,t)|
\right]^{\frac{\gamma}{2}}dx\right]^{\frac{2}{\gamma}}\leq CN^{-1}
 \left[\int_{|x|> \frac{1}{N}}|x|^{-\gamma}dx\right]^{\frac{2}{\gamma}}
 \leq C N^{-1}N^{2-\frac{2}{\gamma}}=N^{\frac{\gamma-2}{\gamma}}.
\end{align}
 We define
$-i\beta N^{3}\eta^{3}t+ixN\eta-\frac{i\gamma t}{N\eta }=-it\beta
 N^{3}(\eta^{3}+\frac{\gamma}{\beta \eta N^{4}}-\frac{x\eta}{\beta N^{2}t}):
=itN^{3}\beta\phi_{6}(\eta)$.Obviously,
  we have   $|\phi_{6}^{(3)}(\eta)|\geq 1$.  By using  Lemma 2.7, we have
\begin{align}\label{2.040}
 |K(x,t)|=N\left|\int_{1}^{4}e^{-iN^{3}\eta^{3}t-\frac{i\gamma t}{N\eta }+ixN\eta}d\eta\right|
 \leq CN(N^{3}t)^{-\frac{1}{3}}
 \leq \frac{C}{t^{\frac{1}{3}}}\leq C\frac{N^{\frac{2}{3}}}{|x|^{\frac{1}{3}}}.
 \end{align}
  Thus, by using  (\ref{2.041}),  for  $\gamma \geq7,$  we have
 \begin{align}\label{2.041}
\left[\int_{|x|>\frac{1}{N}}\left[\sup\limits_{t\in \Omega_{x,3}}
|K(x,t)|\right]^{\frac{\gamma}{2}}dx\right]^{\frac{2}{\gamma}}\leq C
 \left[\int_{|x|>\frac{1}{N}}N^{\frac{\gamma}{3}}
 |x|^{-\frac{\gamma}{6}}dx\right]^{\frac{2}{\gamma}}\leq C
 N^{\frac{\gamma-2}{\gamma}}.
\end{align}
Putting together  (\ref{2.038}), (\ref{2.039}) with (\ref{2.041}),  we derive
\begin{align*}
\|K\|_{L_{x}^{\frac{\gamma}{2}}L_{t}^{\infty}}\leq C N^{\frac{\gamma-2}{\gamma}}.
\end{align*}

We have completed the proof of  Lemma 2.8.

\noindent {\bf Remark3:}Inspired by  the idea of Proposition 2.5 of \cite{Porn},
 we show Lemma 2.8.

\begin{Lemma}\label{lem2.9}
 We assume that    $\gamma\geq 4$ and  $supp(|\mathscr{F}_{x}{f}|)
 \subseteq [N,4N]$,  where $N\in 2^{Z}$, $N\geq a$.  For  $f\in L_{x}^{2}(\R)$,   we have
\begin{align}\label{2.042}
\|U^{\gamma,\beta}(t)f(x)\|_{L_{x}^{\gamma}L_{t}^{\infty}}\leq C N^{\frac{\gamma-2}{2\gamma}}
\|f\|_{L_{x}^{2}}.
\end{align}
Here
\begin{eqnarray}
U^{\gamma,\beta}(t)f(x)=\int_{\SR}e^{ix\xi-it\phi(\xi)}\mathscr{F}_{x}f(\xi)\chi_{[N,4N]}(|\xi|)d\xi.\label{2.043}
\end{eqnarray}
We  assume that  $b>\frac{1}{2},$  $\gamma\geq 4$ and  $supp(|\mathscr{F}_{x}{u}|)
\subseteq [N,4N]$,
where $N\in 2^{Z}$,  $N\geq a.$  Then,   we derive
\begin{eqnarray}
&&\|u\|_{L_{x}^{\gamma}L_{t}^{\infty}}\leq C N^{\frac{\gamma-2}{2\gamma}}
\|u\|_{X_{0,b}}.\label{2.044}
\end{eqnarray}
We  assume  that  $b>\frac{1}{2},$  $\gamma\geq 4$ and  $supp(|\mathscr{F}_{x}{u}|)
\subseteq [a,+\infty)$ for any $t$. Then,  we have
\begin{eqnarray}
&&\|u\|_{L_{x}^{\gamma}L_{t}^{\infty}}\leq C
\|u\|_{X_{s,b}}.\label{2.045}
\end{eqnarray}
Here $s=\frac{1}{2}-\frac{1}{\gamma}+\epsilon.$

\end{Lemma}
\noindent{\bf Proof.}From (2.1) of \cite{GH}, we have
\begin{eqnarray}
\|U^{\gamma,\beta}(t)f(x)\|_{L_{x}^{4}L_{t}^{\infty}}\leq C N^{\frac{1}{4}}
\|f\|_{L_{x}^{2}}\label{2.046}.
\end{eqnarray}
We firstly show that (\ref{2.042}) is valid for  $\gamma \geq 7.$
  We define $Tf:=U^{\gamma,\beta}(t)f$, where $T:L_{x}^{2}\rightarrow L_{x}^{\gamma}L_{t}^{\infty}$.
  Obviously, we have
 $T^{*}F=\int_{\SR}e^{t(\beta\partial_{x}^{3}+\gamma\partial_{x}^{-1})}F dt$.
 By Using the $TT^{*}$ idea, we know  that
 that (\ref{2.042}) is equivalent to
 \begin{align}\label{2.047}
 \left\|\int_{\SR}U^{\gamma,\beta}(t-s)Fds\right\|_{L_{x}^{\gamma}L_{t}^{\infty}}\leq C
  N^{\frac{\gamma-2}{\gamma}}\|F\|_{L_{x}^{\frac{\gamma}{\gamma-1}}L_{t}^{1}}.
 \end{align}
 Here, $F\in L_{x}^{2}L_{t}^{1}(\R\times\R)$ possesses the same frequency support as $u$.
 Thus, to obtain (\ref{2.042}),  it suffices to prove (\ref{2.047}).
 By using a  direct computation, we have
 \begin{align*}
 \mathscr{F}_{x}^{-1}\left(e^{-i(t-s)\phi(\xi)}\mathscr{F}_{x}{F}(\xi,s)\right)
 &=C\int_{\SR}e^{-i(t-s)\phi(\xi)+ix\xi}\mathscr{F}_{x}{F}(\xi,s)d\xi\\
 &=\mathscr{F}_{x}^{-1}\left(e^{-i(t-s)\phi(\xi)}\chi_{[N,4N]}\left(\xi\right)\right)*F(x,s).
 \end{align*}
Thus,  the term on the left hand side of (\ref{2.047}) can be rewritten as
 \begin{align*}
 &\int_{\SR}\mathscr{F}_{x}^{-1}\left(e^{-i(t-s)\phi(\xi)}\chi_{[N,4N]}\left(|\xi|\right)\right)*F(x,s)ds\nonumber\\
 &=\mathscr{F}_{x}^{-1}\left(e^{-it\phi(\xi)}\chi_{[N,4N]}\left(|\xi|\right)\right)*F(x,t)=CK*F.
 \end{align*}
 where $*$ denotes the convolution with respect to variables $x,t$  and
\begin{align}\label{2.048}
 K(x,t)=\int_{\SR}e^{-it\phi(\xi)+ix\xi}\chi_{[N,4N]}\left(|\xi|\right)d\xi.
 \end{align}
From Young inequality and Lemma 2.8, we derive
 \begin{align*}
\|K*F\|_{L_{x}^{\gamma}L_{t}^{\infty}}\leq \|K\|_{L_{x}^{\frac{\gamma}{2}}L_{t}^{\infty}}
\|F\|_{L_{x}^{\frac{\gamma}{\gamma-1}}L_{t}^{1}}\leq C
  N^{\frac{\gamma-2}{\gamma}}\|F\|_{L_{x}^{\frac{\gamma}{\gamma-1}}L_{t}^{1}}.
 \end{align*}
Consequently, when $\gamma\geq7,$ we derive
 \begin{align}\label{2.049}
\|U^{\gamma,\beta}(t)f(x)\|_{L_{x}^{\gamma}L_{t}^{\infty}}\leq C N^{\frac{\gamma-2}{2\gamma}}
\|f\|_{L_{x}^{2}}.
\end{align}
Interpolating (\ref{2.046}) with (\ref{2.049}) leads to
\begin{align}\label{2.050}
\|U^{\gamma,\beta}(t)f(x)\|_{L_{x}^{\gamma}L_{t}^{\infty}}\leq C N^{\frac{\gamma-2}{2\gamma}}
\|f\|_{L_{x}^{2}},4\leq \gamma \leq 7.
\end{align}
By  changing variable $\tau=\lambda-\phi(\xi)$, we derive
 \begin{eqnarray}
&&u(x,t)=\frac{1}{2\pi}\int_{\SR^{2}}e^{ix\xi+it\tau}\mathscr{F}u(\xi,\tau)
d\xi d\tau\nonumber\\
&&=\frac{1}{2\pi}\int_{\SR^{2}}e^{ix\xi+it(\lambda-\phi(\xi))}
\mathscr{F}u(\xi,\lambda-\phi(\xi))d\xi d\lambda\nonumber\\
&&=\frac{1}{2\pi}\int_{\SR}e^{it\lambda}\left(\int_{\SR}e^{ix\xi-it\phi(\xi)}
\mathscr{F}u(\xi,\lambda-\phi(\xi)) d\xi\right)d\lambda\label{2.051}.
 \end{eqnarray}
By using  (\ref{2.050}),  (\ref{2.051}) and Minkowski's inequality
 and the change of variable $\lambda=\tau+\phi,$  for $b>\frac{1}{2}$,     we have
\begin{eqnarray}
&&\|u\|_{L_{x}^{\gamma}L_{t}^{\infty}}\leq C\int_{\SR}\left\|\left(\int_{\SR}e^{ix\xi-it\phi(\xi)}
\mathscr{F}u(\xi,\lambda-\phi(\xi))d\xi\right)
\right\|_{L_{x}^{\gamma}L_{t}^{\infty}} d\lambda\nonumber\\&&\leq CN^{\frac{\gamma-2}{2\gamma}}\int_{\SR}\left\|\mathscr{F}u(\xi,\lambda-\phi(\xi))\right\|_{L^{2}}
d\lambda\nonumber\\&&\leq CN^{\frac{\gamma-2}{2\gamma}}
\left[\int_{\SR}(1+|\lambda|)^{2b}\left\|\mathscr{F}u(\xi,\lambda-\phi(\xi))\right\|_{L^{2}}^{2}
d\lambda\right]^{\frac{1}{2}}
\left[\int_{\SR}(1+|\lambda|)^{-2b}d\lambda\right]^{\frac{1}{2}}\nonumber\\
&&\leq CN^{\frac{\gamma-2}{2\gamma}}\left[\int_{\SR}(1+|\tau+\phi(\xi)|)^{2b}\left\|\mathscr{F}u(\xi,\tau)
\right\|_{L^{2}}^{2}d\tau\right]^{\frac{1}{2}}
=CN^{\frac{\gamma-2}{2\gamma}}\|u\|_{X_{0,b}}.\label{2.052}
\end{eqnarray}
Since
\begin{eqnarray}
P_{a}U^{\gamma,\beta}u_{0}=\sum_{N\geq a}\int_{N}^{4N}e^{ix\xi-it\phi}\mathscr{F}_{x}u_{0}d\xi,\label{2.053}
\end{eqnarray}
thus, by using  the Minkowski  equality  and  Cauchy-Schwarz inequality,   we have
\begin{eqnarray}
&&\left\|P_{a}U^{\gamma,\beta}u_{0}\right\|_{L_{x}^{\gamma}L_{t}^{\infty}}\leq \sum_{N\geq a}\left\|\int_{N}^{4N}e^{ix\xi-it\phi}\mathscr{F}_{x}u_{0}d\xi\right\|_{L_{x}^{\gamma}L_{t}^{\infty}}\leq C\sum_{N\geq a}N^{\frac{1}{2}-\frac{1}{\gamma}}\|\chi_{[N,4N]}(|\xi|)\mathscr{F}_{x}u_{0}\|_{L_{\xi}^{2}}\nonumber\\&&\leq C\sum_{N\geq a}\left[N^{2(\frac{1}{2}-\frac{1}{\gamma}+\epsilon)}\|\chi_{[N,4N]}(|\xi|)\mathscr{F}_{x}u_{0}\|_{L_{\xi}^{2}}^{2}\right]^{\frac{1}{2}}\nonumber\\
&&\leq C\|u_{0}\|_{H^{s}},\label{2.054}
\end{eqnarray}
where $s=\frac{1}{2}-\frac{1}{\gamma}+\epsilon.$
By using (\ref{2.054}) and a proof similar to (\ref{2.044}), we obtain that (\ref{2.045}) is valid.

We have completed the proof of  Lemma 2.9.

\begin{Lemma}\label{lem2.10}
We assume  that     $\gamma\geq 4$ and
  $supp(|\mathscr{F}_{x}{u}|)\subseteq [N,4N]$, $N\in 2^{Z}$.  For $b>\frac{1}{2}$ and $f\in L_{x}^{2}(\R)$,    we have
\begin{align}\label{2.055}
\|\psi (t)u\|_{L_{xt}^{\infty}}\leq C N^{\frac{1}{4}-\epsilon}\|u\|_{X_{0,b}}.
\end{align}
For $b>\frac{1}{2},$   and  $supp(|\mathscr{F}_{x}{u}|)
\subseteq (0,a]$ for any $t$, we have
\begin{eqnarray}
&&\|\psi(t)u\|_{L_{xt}^{\infty}}\leq C
\|D_{x}^{s}u\|_{X_{0,b}}.\label{2.056}
\end{eqnarray}
Here $s=\frac{1}{4}-\epsilon.$  For $b>\frac{1}{2},$   and  $supp(|\mathscr{F}_{x}{u}|)
\subseteq (0,a]$ for any $t$, we have
\begin{eqnarray}
&&\|\psi(t)u\|_{L_{x}^{\frac{2}{1-2\epsilon}}L_{t}^{\infty}}\leq C
\|D_{x}^{-\frac{1}{4}}u\|_{X_{0,b}}.\label{2.057}
\end{eqnarray}
\end{Lemma}
\noindent{\bf Proof.} Without loss of generality,  we can assume that $supp(\mathscr{F}_{x}{u})\subseteq [N,4N]$.
 By using
 the Cauchy-Schwarz inequality and the Minkowski's inequality, we have
\begin{eqnarray}
&&\left|u(x)\right|=\left|\int_{N}^{4N}e^{ix\xi}\mathscr{F}_{x}u(\xi,t)d\xi\right|
\leq \int_{N}^{4N}\left|\mathscr{F}_{x}u(\xi,t)\right|d\xi\leq
CN^{\frac{1}{2}}\left[\int_{N}^{4N}\left|\mathscr{F}_{x}u(\xi,t)
\right|^{2}d\xi\right]^{\frac{1}{2}}\nonumber\\
&&=CN^{\frac{1}{2}}\|\mathscr{F}_{x}u(\xi,t)\|_{L_{\xi}^{2}}=
CN^{\frac{1}{2}}\|u\|_{L_{x}^{2}}\label{2.058}.
\end{eqnarray}
Thus, by  using  the  Minkowski's  inequality,   from (\ref{2.058}) and (\ref{2.017}), we have
\begin{eqnarray}
&&\left\|\psi(t)u\right\|_{L_{xt}^{\infty}}\leq CN^{\frac{1}{2}}
\|u\|_{L_{t}^{\infty}L_{x}^{2}}\leq
 CN^{\frac{1}{2}}\|u\|_{L_{x}^{2}L_{t}^{\infty}}\leq C
 N^{\frac{1}{4}-\epsilon}\|u\|_{X_{0,b}}\label{2.059}.
\end{eqnarray}
By using (\ref{2.059}) and a proof similar to (\ref{2.045}), we obtain that (\ref{2.056}) is valid.
Interpolating (\ref{2.09}) with (\ref{2.056}) yields (\ref{2.057}).

This completes the proof of Lemma 2.10.
\begin{Lemma}\label{lem2.11}
For $b>\frac{1}{2}$, we have
\begin{eqnarray}
&&\left\|D_{x}^{-\frac{1}{2}-4\epsilon}\mathscr{F}_{x}^{-1}
\left(\chi_{|\xi|\geq a}(\xi)\mathscr{F}u(\xi)\right)\right\|_{L_{xt}^{\infty}}\leq C\left\|u\right\|_{X_{0,b}}.\label{2.060}
\end{eqnarray}
\end{Lemma}
\noindent{\bf Proof.}
Using the Sobolev embeddings theorem and (\ref{2.03}),we have
\begin{eqnarray}
&&\left\|D_{x}^{-\frac{1}{2}-4\epsilon}\mathscr{F}_{x}^{-1}\left(\chi_{|\xi|\geq a}(\xi)\mathscr{F}U^{\gamma,\beta}(t)u_{0}\right)\right\|_{L_{xt}^{\infty}}\nonumber\\
&&\leq C\left\|D_{x}^{-\frac{1}{2}-4\epsilon}J_{x}^{\frac{1}{8}+\epsilon}J_{t}^{\frac{1}{8}+\epsilon}\mathscr{F}_{x}^{-1}\left(\chi_{|\xi|\geq a}(\xi)\mathscr{F}U^{\gamma,\beta}(t)u_{0}\right)\right\|_{L_{xt}^{8}}\nonumber\\
&&=C\left\|\int_{\SR}e^{ix\xi-it\phi(\xi)}\chi_{|\xi|\geq a}(\xi)\left|\xi\right|^{-\frac{1}{2}-4\epsilon}\langle\xi\rangle^{\frac{1}{8}+\epsilon}\left(1+\left|\beta\xi^{3}
+\frac{\gamma}{\xi}\right|\right)^{\frac{1}{8}+\epsilon}\mathscr{F}_{x}u_{0}d\xi\right\|_{L_{xt}^{8}}\nonumber\\
&&\leq C\left\|\left|\xi\right|^{-\frac{1}{2}-4\epsilon}\langle\xi\rangle^{\frac{1}{8}+\epsilon}\left|\beta\xi^{3}
+\frac{\gamma}{\xi}\right|^{\frac{1}{8}+\epsilon}\chi_{|\xi|\geq a}(\xi)\mathscr{F}_{x}u_{0}\right\|_{L_{\xi}^{2}}\nonumber\\
&&\leq C\left\|\mathscr{F}_{x}u_{0}\right\|_{L_{\xi}^{2}}=C\left\|u_{0}\right\|_{L_{x}^{2}}.\label{2.061}
\end{eqnarray}
By using (\ref{2.061}) and a proof similar to (\ref{2.04}), we obtain that (\ref{2.060}) is valid.

This completes the proof of Lemma 2.11.
\begin{Lemma}\label{lem2.12}
Let $s>0$, $1<p<\infty$.  Then, we have
\begin{eqnarray}
&&\left\|J^{s}(fg)-fJ^{s}g\right\|_{L^{p}}\leq C\left\|f_{x}\right\|_{L^{\infty}}\left\|J^{s-1}g\right\|_{L^{p}}
+C\left\|J^{s}f\right\|_{L^{p}}\left\|g\right\|_{L^{\infty}}.\label{2.062}
\end{eqnarray}
\end{Lemma}

For the  proof  of  Lemma 2.12, we refer the readers to  Lemma XI of \cite{KP}.

\bigskip
\bigskip

\noindent{\large\bf 3. Multilinear estimates}

\setcounter{equation}{0}

 \setcounter{Theorem}{0}

\setcounter{Lemma}{0}

 \setcounter{section}{3}
 In this section, we present some crucial multilinear estimates which play an
  important role in establishing Theorems 1.1 and Theorem 1.2.

\begin{Lemma}\label{Lemma3.1}
Let $s\geq\frac{1}{2}-\frac{2}{k}+2\epsilon,k\geq5$
and $b=\frac{1}{2}+\frac{\epsilon}{24}$ and $g_{j}=\psi(t)u_{j}.$
Then, we have
\begin{eqnarray}
\left\|\partial_{x}(\prod\limits_{j=1}^{k+1}g_{j})
\right\|_{X_{s,-\frac{1}{2}+\frac{\epsilon}{12}}}
\leq C\prod_{j=1}^{k+1}\|g_{j}\|_{X_{s,b}}.\label{3.01}
\end{eqnarray}
\end{Lemma}
\noindent {\bf Proof.}To prove (\ref{3.01}),
by duality,  it suffices to prove
\begin{eqnarray}
\left|\int_{\SR^{2}}J^{s}\partial_{x}
\left(\prod_{j=1}^{k+1}g_{j}\right)\bar{h}dxdt
\right|\leq C\|h\|_{X_{0,\frac{1}{2}-\frac{\epsilon}{12}}}
\prod_{j=1}^{k+1}\|g_{j}\|_{X_{s,b}}\label{3.02}.
\end{eqnarray}
We define
\begin{eqnarray*}
f(\xi,\tau)=\langle\sigma\rangle^{\frac{1}{2}-\frac{\epsilon}{12}} \mathscr{F}h(\xi,\tau),
f_{j}(\xi_{j},\tau_{j})=\langle\xi_{j}\rangle^{s}\langle\sigma_{j}\rangle^{b}
\mathscr{F}g_{j}(\xi_{j},\tau_{j})(1\leq j\leq k+1).
\end{eqnarray*}
By using the Plancherel identity, to prove (\ref{3.02}), it suffices to prove
\begin{eqnarray}
\int_{\xi=\sum\limits_{j=1}^{k+1}\xi_{j}}\int_{\tau=\sum\limits_{j=1}^{k+1}\tau_{j}}
\frac{|\xi|\langle\xi\rangle^{s}f\prod\limits_{j=1}^{k+1}f_{j}}
{\langle\sigma\rangle^{\frac{1}{2}-\frac{\epsilon}{12}}\prod\limits_{j=1}^{k+1}
\langle\xi_{j}\rangle^{s}\langle\sigma_{j}\rangle^{b}}d\delta
&&\leq C\|f\|_{L^{2}}\prod\limits_{j=1}^{k+1}\|f_{j}\|_{L^{2}},\label{3.03}
\end{eqnarray}
where $d\delta=d\xi_{1}d\xi_{2}\cdot \cdot \cdot d\xi_{k}
d\xi d\tau_{1}d\tau_{2}\cdot \cdot \cdot d\tau_{k}d\tau.$

\noindent We define
\begin{eqnarray*}
&&K(\xi_{1},\xi_{2},\cdot\cdot\cdot,\xi_{k},\xi,\tau_{1},\tau_{2},\cdot\cdot\cdot,\tau_{k},\tau)
=\frac{|\xi|\langle\xi\rangle^{s}}
{\langle\sigma\rangle^{\frac{1}{2}-\frac{\epsilon}{12}}\prod\limits_{j=1}^{k+1}
\langle\xi_{j}\rangle^{s}\langle\sigma_{j}\rangle^{b}},\nonumber\\&&
\mathscr{F}F=\frac{f}{\langle\sigma\rangle^{\frac{1}{2}-\frac{\epsilon}{12}}},
\mathscr{F}F_{j}=\frac{f_{j}}{\langle\sigma_{j}\rangle^{b}}
(1\leq j\leq k+1),\\
&&I_{1}=\int_{\xi=\sum\limits_{j=1}^{k+1}\xi_{j}}
\int_{\tau=\sum\limits_{j=1}^{k+1}\tau_{j}}K(\xi_{1},\xi_{2},
\cdot\cdot\cdot,\xi_{k},\xi,\tau_{1},\tau_{2},\cdot\cdot\cdot,\tau_{k},\tau)
f\prod\limits_{j=1}^{k+1}f_{j}d\delta.
\end{eqnarray*}
Without loss of generality, we can assume that $|\xi_{1}|\geq
 |\xi_{2}|\geq\cdot\cdot\cdot\geq |\xi_{k+1}|.$

\noindent Obviously,
\begin{eqnarray*}
\Omega:=\left\{(\xi_{1},\xi_{2},\cdot\cdot\cdot,\xi_{k},\xi,\tau_{1},\tau_{2},
\cdot\cdot\cdot,\tau_{k},\tau)\in \R^{2(k+1)}:|\xi_{1}|\geq |\xi_{2}|\geq
\cdot\cdot\cdot\geq |\xi_{k+1}|\right\}\subset \bigcup\limits_{j=0}^{k+1}\Omega_{j}.
\end{eqnarray*}
Here,
\begin{eqnarray*}
&&\Omega_{0}=\left\{(\xi_{1},\xi_{2},\cdot\cdot\cdot,\xi_{k},\xi,\tau_{1},\tau_{2},
\cdot\cdot\cdot,\tau_{k},\tau)\in \Omega,|\xi_{1}|\leq80ka\right\},\\
&&\Omega_{1}=\left\{(\xi_{1},\xi_{2},\cdot\cdot\cdot,\xi_{k},\xi,\tau_{1},\tau_{2},
\cdot\cdot\cdot,\tau_{k},\tau)\in \Omega,|\xi_{1}|\geq80ka,|\xi_{1}|\geq80ka |\xi_{2}|\right\},\\
&&\Omega_{2}=\left\{(\xi_{1},\xi_{2},\cdot\cdot\cdot,\xi_{k},\xi,\tau_{1},\tau_{2},
\cdot\cdot\cdot,\tau_{k},\tau)\in \Omega,|\xi_{1}|\geq80ka, |\xi_{1}|\sim |\xi_{2}|\geq 80ka|\xi_{3}|\right\},\\
&&\Omega_{3}=\left\{(\xi_{1},\xi_{2},\cdot\cdot\cdot,\xi_{k},\xi,\tau_{1},\tau_{2},
\cdot\cdot\cdot,\tau_{k},\tau)\in \Omega,|\xi_{1}|\geq80ka, |\xi_{1}|\sim
 |\xi_{3}|\geq 80ka|\xi_{4}|\right\},\\
&&\Omega_{4}=\left\{(\xi_{1},\xi_{2},\cdot\cdot\cdot,\xi_{k},\xi,\tau_{1},\tau_{2},
\cdot\cdot\cdot,\tau_{k},\tau)\in \Omega,|\xi_{1}|\geq80ka, |\xi_{1}|\sim|\xi_{4}|\geq 80ka|\xi_{5}|\right\},\\
&&\Omega_{5}=\left\{(\xi_{1},\xi_{2},\cdot\cdot\cdot,\xi_{k},\xi,\tau_{1},\tau_{2},
\cdot\cdot\cdot,\tau_{k},\tau)\in \Omega,|\xi_{1}|\geq80ka, |\xi_{1}|\sim  |\xi_{5}|\geq 80ka|\xi_{6}|\right\},\\
&&\Omega_{6}=\left\{(\xi_{1},\xi_{2},\cdot\cdot\cdot,\xi_{k},\xi,\tau_{1},\tau_{2},
\cdot\cdot\cdot,\tau_{k},\tau)\in \Omega,|\xi_{1}|\geq80ka, |\xi_{1}|\sim |\xi_{l-1}|\geq 80ka|\xi_{l}|(7\leq l\leq k+1)\right\},\\
&&\Omega_{7}=\left\{(\xi_{1},\xi_{2},\cdot\cdot\cdot,\xi_{k},\xi,\tau_{1},\tau_{2},
\cdot\cdot\cdot,\tau_{k},\tau)\in \Omega,|\xi_{1}|\geq80ka, |\xi_{1}|\sim  |\xi_{k+1}|\right\}.
\end{eqnarray*}
(1)Case $(\xi_{1},\xi_{2},\cdot\cdot\cdot,\xi_{k},\xi,\tau_{1},\tau_{2},\cdot\cdot\cdot,
\tau_{k},\tau)\in \Omega_{0}$,
 by using the Plancherel identity,  the H\"older inequality and  (\ref{2.04}), we have
\begin{eqnarray}
&&I_{1}\leq C\int_{\xi=\sum\limits_{j=1}^{k+1}\xi_{j}}\int_{\tau=\sum\limits_{j=1}^{k+1}\tau_{j}}
\frac{f\prod\limits_{j=1}^{k+1}f_{j}}
{\langle\sigma\rangle^{\frac{1}{2}-\frac{\epsilon}{12}}\prod\limits_{j=1}^{k+1}
\langle\xi_{j}\rangle^{s}\langle\sigma_{j}\rangle^{b}}d\delta\nonumber\\
&&\leq C\|F\|_{L_{xt}^{2}}\prod\limits_{j=1}^{k+1}\|F_{j}\|_{L_{xt}^{2(k+1)}}
\leq C\|f\|_{L^{2}}\prod\limits_{j=1}^{k+1}\|F_{j}\|_{X_{\frac{1}{2}-\frac{2}{k+1},b}}
\leq C\|f\|_{L^{2}}\prod\limits_{j=1}^{k+1}\|f_{j}\|_{L^{2}}\label{3.04}.
\end{eqnarray}
(2)Case $(\xi_{1},\xi_{2},\cdot\cdot\cdot,\xi_{k},\xi,\tau_{1},\tau_{2},\cdot\cdot\cdot,\tau_{k},\tau)\in \Omega_{1}$, we consider
\begin{eqnarray}
&&|1+\frac{\gamma}{3\beta\xi_{1}^{2}\xi_{2}^{2}}|\geq\frac{1}{2},\label{3.05}\\
&&|1+\frac{\gamma}{3\beta\xi_{1}^{2}\xi_{2}^{2}}|<\frac{1}{2}.\label{3.06}
\end{eqnarray}
When (\ref{3.05}) is valid, we have
\begin{eqnarray}
&&K(\xi_{1},\xi_{2},\cdot\cdot\cdot,\xi_{k},\xi,\tau_{1},\tau_{2},\cdot\cdot\cdot,\tau_{k},\tau)\leq C
\frac{|\xi_{1}^{2}-\xi_{2}^{2}|^{\frac{1}{2}}}
{\langle\sigma\rangle^{\frac{1}{2}-\frac{\epsilon}{12}}\prod\limits_{j=2}^{k+1}
\langle\xi_{j}\rangle^{s}\prod\limits_{j=1}^{k+1}\langle\sigma_{j}\rangle^{b}}\nonumber\\
&&\leq C \frac{|\phi^{\prime}(\xi_{1})-\phi^{\prime}(\xi_{2})|^{\frac{1}{2}}}
{\langle\sigma\rangle^{\frac{1}{2}-\frac{\epsilon}{12}}\prod\limits_{j=2}^{k+1}
\langle\xi_{j}\rangle^{s}\prod\limits_{j=1}^{k+1}\langle\sigma_{j}\rangle^{b}}
.\label{3.07}
\end{eqnarray}
By using  (\ref{3.07}), the Plancherel identity,  the H\"older inequality,
Lemma 2.5, (\ref{2.04}) and (\ref{2.07}),
we have
\begin{eqnarray*}
&&I_{1}\leq C\int_{\xi=\sum\limits_{j=1}^{k+1}\xi_{j}}\int_{\tau=\sum\limits_{j=1}^{k+1}\tau_{j}}
\frac{f(\prod\limits_{j=1}^{k+1}f_{j})|\phi^{\prime}(\xi_{1})-\phi^{\prime}(\xi_{2})|^{\frac{1}{2}}\prod\limits_{j=2}^{k+1}
\langle\xi_{j}\rangle^{-s}}{\langle\sigma\rangle^{\frac{1}{2}-\frac{\epsilon}{12}}
\prod\limits_{j=1}^{k+1}\langle\sigma_{j}\rangle^{b}}d\delta\nonumber\\
&&\leq C
\left\|I^{1/2}(F_{1},F_{2})\right\|_{L_{xt}^{2}}
\left(\prod_{j=3}^{k+1}\|J_{x}^{-\frac{ks}{k-1}}F_{j}\|_{L_{xt}^{\frac{8(k-1)}{3-\epsilon}}}\right)
\|F\|_{L_{xt}^{\frac{8}{1+\epsilon}}}\nonumber\\
&&\leq C\|F\|_{X_{0,\frac{1}{2}-\frac{\epsilon}{12}}}\prod\limits_{j=1}^{k+1}\|F_{j}\|_{X_{0,b}}\leq
C\left(\prod_{j=1}^{k+1}\|f_{j}\|_{L_{\xi\tau}^{2}}\right)\|f\|_{L_{\xi\tau}^{2}}.
\end{eqnarray*}
When (\ref{3.06})  is valid,  we have $|\xi_{1}|\sim|\xi_{2}|^{-1}$. In this case,  we consider
\begin{eqnarray}
&&\left|1+\frac{\gamma}{3\beta\xi^{2}\xi_{k+1}^{2}}\right|\geq\frac{1}{2},\label{3.08}\\
&&\left|1+\frac{\gamma}{3\beta\xi^{2}\xi_{k+1}^{2}}\right|<\frac{1}{2},\label{3.09}
\end{eqnarray}
respectively.
When (\ref{3.08})  is valid,    we have
\begin{eqnarray}
&&K(\xi_{1},\xi_{2},\cdot\cdot\cdot,\xi_{k},\xi,\tau_{1},\tau_{2},
\cdot\cdot\cdot,\tau_{k},\tau)\leq C
\frac{|\xi_{1}||\xi_{2}|^{\frac{1}{4}+\epsilon}|\xi|^{1-4\epsilon}}
{\langle\sigma\rangle^{\frac{1}{2}-\frac{\epsilon}{12}}
\prod\limits_{j=1}^{k+1}\langle\sigma_{j}\rangle^{b}}
.\label{3.010}
\end{eqnarray}
By using  (\ref{3.010}), the Plancherel identity,  the H\"older inequality,
Lemma 2.5, (\ref{2.08}),  (\ref{2.09})   and   (\ref{2.056}),
we have
\begin{eqnarray*}
&&I_{1}\leq C\int_{\xi=\sum\limits_{j=1}^{k+1}\xi_{j}}\int_{\tau=\sum\limits_{j=1}^{k+1}\tau_{j}}
\frac{f(\prod\limits_{j=1}^{k+1}f_{j})|\phi^{\prime}(\xi)-\phi^{\prime}
(\xi_{k+1})|^{\frac{1}{2}-2\epsilon}|\xi_{1}||\xi_{2}|^{\frac{1}{4}+\epsilon}}
{\langle\sigma\rangle^{\frac{1}{2}-\frac{\epsilon}{12}}
\prod\limits_{j=1}^{k+1}\langle\sigma_{j}\rangle^{b}}d\delta\nonumber\\&&\leq C
\left\|I^{\frac{1}{2}-2\epsilon}(F_{k+1},F)\right\|_{L_{xt}^{2}}
\left(\prod_{j=3}^{k}\|F_{j}\|_{L_{xt}^{\infty}}\right)
\|D_{x}F_{1}\|_{L_{x}^{\infty}L_{t}^{2}}\|D_{x}^{\frac{1}{4}+\epsilon}
F_{2}\|_{L_{x}^{2}L_{t}^{\infty}}\nonumber\\
&&\leq C\|F\|_{X_{0,\frac{1}{2}-\frac{\epsilon}{12}}}
\prod\limits_{j=1}^{k+1}\|F_{j}\|_{X_{0,b}}\leq
C\left(\prod_{j=1}^{k+1}\|f_{j}\|_{L_{\xi\tau}^{2}}\right)\|f\|_{L_{\xi\tau}^{2}}.
\end{eqnarray*}
When (\ref{3.09})  is valid,  we have that  $|\xi|\sim|\xi_{k+1}|^{-1}$, thus,  we have
\begin{eqnarray}
&&K(\xi_{1},\xi_{2},\cdot\cdot\cdot,\xi_{k},\xi,\tau_{1},\tau_{2},\cdot\cdot\cdot,\tau_{k},\tau)\leq C
\frac{|\xi_{1}||\xi_{2}|^{\frac{1}{4}+\epsilon}|\xi_{k+1}|^{\frac{1}{4}}|\xi|^{1-2\epsilon}}
{\langle\sigma\rangle^{\frac{1}{2}-\frac{\epsilon}{12}}\prod\limits_{j=1}^{k+1}\langle\sigma_{j}\rangle^{b}}
.\label{3.011}
\end{eqnarray}
By using  (\ref{3.011}), the Plancherel identity,  the H\"older inequality, (\ref{2.08})-(\ref{2.010}), (\ref{2.056}), (\ref{2.057}),
we have
\begin{eqnarray*}
&&I_{1}\leq C\int_{\xi=\sum\limits_{j=1}^{k+1}\xi_{j}}\int_{\tau=\sum\limits_{j=1}^{k+1}\tau_{j}}
\frac{f(\prod\limits_{j=1}^{k+1}f_{j})|\xi_{1}||\xi_{2}|^{\frac{1}{4}}|\xi_{k+1}|^{\frac{1}{4}}|\xi|^{1-2\epsilon}}
{\langle\sigma\rangle^{\frac{1}{2}-\frac{\epsilon}{12}}
\prod\limits_{j=1}^{k+1}\langle\sigma_{j}\rangle^{b}}d\delta\nonumber\\
&&\leq C
\|D_{x}F_{1}\|_{L_{x}^{\infty}L_{t}^{2}}\|D_{x}^{\frac{1}{4}+\epsilon}F_{2}\|_{L_{x}^{2}L_{t}^{\infty}}
\left(\prod_{j=3}^{k}\|F_{j}\|_{L_{xt}^{\infty}}\right)\|D_{x}^{\frac{1}{4}}F_{k+1}\|_{L_{x}^{\frac{2}{1-2\epsilon}}L_{t}^{\infty}}
\|D_{x}^{1-2\epsilon}F\|_{L_{x}^{\frac{1}{\epsilon}}L_{t}^{2}}\nonumber\\
&&\leq C\|F\|_{X_{0,\frac{1}{2}-\frac{\epsilon}{12}}}\prod\limits_{j=1}^{k+1}\|F_{j}\|_{X_{0,b}}\leq
C\left(\prod_{j=1}^{k+1}\|f_{j}\|_{L_{\xi\tau}^{2}}\right)\|f\|_{L_{\xi\tau}^{2}}.
\end{eqnarray*}
(3)Case $(\xi_{1},\xi_{2},\cdot\cdot\cdot,\xi_{k},\xi,\tau_{1},\tau_{2},
\cdot\cdot\cdot,\tau_{k},\tau)\in \Omega_{2}$,  we consider
\begin{eqnarray}
&&|1+\frac{\gamma}{3\beta\xi_{2}^{2}\xi_{3}^{2}}|\geq\frac{1}{2},\label{3.012}\\
&&|1+\frac{\gamma}{3\beta\xi_{2}^{2}\xi_{3}^{2}}|<\frac{1}{2}.\label{3.013}
\end{eqnarray}
When (\ref{3.012}) is valid,  since $s\geq\frac{1}{2}-\frac{2}{k}+2\epsilon$,  we have
\begin{align}
&K(\xi_{1},\xi_{2},\cdot\cdot\cdot,\xi_{k},\xi,\tau_{1},\tau_{2},\cdot\cdot\cdot,\tau_{k},\tau)\leq
C\frac{|\xi_{2}^{2}-\xi_{3}^{2}|^{\frac{1}{2}}}
{\langle\sigma\rangle^{\frac{1}{2}-\frac{\epsilon}{12}}\prod\limits_{j=2}^{k+1}
\langle\xi_{j}\rangle^{s}\prod\limits_{j=1}^{k+1}\langle\sigma_{j}\rangle^{b}}
\notag\\
&\leq C\frac{|\xi_{2}^{2}-\xi_{3}^{2}|^{\frac{1}{2}}
|\xi_{1}|^{\frac{1}{6}}\prod\limits_{j=4}^{k+1}
\langle\xi_{j}\rangle^{-\frac{3k-11}{6(k-2)}-2k\epsilon}}{\langle\sigma\rangle^{\frac{1}{2}-\frac{\epsilon}{12}}
\prod\limits_{j=1}^{k+1}\langle\sigma_{j}\rangle^{b}}\nonumber\\
&\leq C\frac{|\phi^{\prime}(\xi_{2})-\phi^{\prime}(\xi_{3})|^{\frac{1}{2}}|\xi_{1}|^{\frac{1}{6}}\prod\limits_{j=4}^{k+1}
\langle\xi_{j}\rangle^{-\frac{3k-11}{6(k-2)}-2k\epsilon}}{\langle\sigma\rangle^{\frac{1}{2}-\frac{\epsilon}{12}}
\prod\limits_{j=1}^{k+1}\langle\sigma_{j}\rangle^{b}}.\label{3.014}
\end{align}
By using (\ref{3.014}), the Plancherel identity, the H\"older inequality,   Lemma 2.5, (\ref{2.04}), (\ref{2.05}),
(\ref{2.07}),  we have
\begin{eqnarray*}
&&I_{1}\leq C\int_{\xi=\sum\limits_{j=1}^{k+1}\xi_{j}}\int_{\tau=\sum\limits_{j=1}^{k+1}\tau_{j}}
\frac{f(\prod\limits_{j=1}^{k+1}f_{j})|\phi^{\prime}(\xi_{2})-\phi^{\prime}(\xi_{3})|^{\frac{1}{2}}|\xi_{1}|^{\frac{1}{6}}
\prod\limits_{j=4}^{k+1}
\langle\xi_{j}\rangle^{-\frac{3k-11}{6(k-2)}-2\epsilon}}{\langle\sigma\rangle^{\frac{1}{2}-\frac{\epsilon}{12}}
\prod\limits_{j=1}^{k+1}\langle\sigma_{j}\rangle^{b}}d\delta\nonumber\\&&\leq C
\left\|I^{1/2}(F_{2},F_{3})\right\|_{L_{xt}^{2}}
\left(\prod_{j=4}^{k+1}\|J_{x}^{-\frac{3k-11}{6(k-2)}-2\epsilon}F_{j}\|_{L_{xt}^{\frac{24(k-2)}{5-3\epsilon}}}\right)
\|D_{x}^{\frac{1}{6}}F_{1}\|_{L_{xt}^{6}}\|F\|_{L_{xt}^{\frac{8}{1+\epsilon}}}\nonumber\\
&&\leq C\left(\prod_{j=1}^{k+1}\|F_{j}\|_{X_{0,b}}\right)
\|F\|_{X_{0,\frac{1}{2}-\frac{\epsilon}{12}}}\leq C\left(\prod_{j=1}^{k+1}\|f_{j}\|_{L_{\xi\tau}^{2}}\right)\|f\|_{L_{\xi\tau}^{2}}.
\end{eqnarray*}
When (\ref{3.013}) is valid,  we have $|\xi_{2}|\sim|\xi_{3}|^{-1}$,
we consider $|\xi|\leq a,|\xi|\geq a$, respectively.

\noindent When $|\xi|\leq a,$  we have
\begin{eqnarray}
&&K(\xi_{1},\xi_{2},\cdot\cdot\cdot,\xi_{k},\xi,\tau_{1},\tau_{2},\cdot\cdot\cdot,\tau_{k},\tau)\leq C
\frac{\prod\limits_{j=1}^{k} \langle \xi_{j}\rangle^{-s}}
{\langle\sigma\rangle^{\frac{1}{2}-\frac{\epsilon}{12}}\prod\limits_{j=1}^{k+1}\langle\sigma_{j}\rangle^{b}}
.\label{3.015}
\end{eqnarray}
By using  (\ref{3.015}), the Plancherel identity,  the H\"older inequality, (\ref{2.04}) and  (\ref{2.056}),
we have
\begin{eqnarray*}
&&I_{1}\leq C\int_{\xi=\sum\limits_{j=1}^{k+1}\xi_{j}}\int_{\tau=\sum\limits_{j=1}^{k+1}\tau_{j}}
\frac{\prod\limits_{j=1}^{k} \langle \xi_{j}\rangle^{-s}f(\prod\limits_{j=1}^{k+1}f_{j})}
{\langle\sigma\rangle^{\frac{1}{2}-\frac{\epsilon}{12}}
\prod\limits_{j=1}^{k+1}\langle\sigma_{j}\rangle^{b}}d\delta\nonumber\\&&\leq C
\left(\prod_{j=1}^{k}\|J_{x}^{-s}F_{j}\|_{L_{xt}^{2k}}\right)\|F_{k+1}\|_{L_{xt}^{\infty}}
\|F\|_{L_{xt}^{2}}\nonumber\\
&&\leq C\|F\|_{X_{0,\frac{1}{2}-\frac{\epsilon}{12}}}\prod\limits_{j=1}^{k+1}\|F_{j}\|_{X_{0,b}}\leq
C\left(\prod_{j=1}^{k+1}\|f_{j}\|_{L_{\xi\tau}^{2}}\right)\|f\|_{L_{\xi\tau}^{2}}.
\end{eqnarray*}
When $|\xi|\geq a,$  we consider
(\ref{3.08}), (\ref{3.09}), respectively.

\noindent When (\ref{3.08}) is valid,  we have
\begin{eqnarray}
&&K(\xi_{1},\xi_{2},\cdot\cdot\cdot,\xi_{k},\xi,\tau_{1},\tau_{2},\cdot\cdot\cdot,\tau_{k},\tau)\leq C
\frac{|\xi_{2}||\xi_{1}|^{-\frac{1}{2}-\epsilon}|\xi_{3}|^{\frac{1}{4}+\epsilon}|\xi|^{1-4\epsilon}}
{\langle\sigma\rangle^{\frac{1}{2}-\frac{\epsilon}{12}}\prod\limits_{j=1}^{k+1}\langle\sigma_{j}\rangle^{b}}
.\label{3.016}
\end{eqnarray}
By using  (\ref{3.016}), the Plancherel identity,  the H\"older inequality,
Lemmas 2.5,  2.11,   (\ref{2.08}),  (\ref{2.09}) and (\ref{2.056}),
we have
\begin{eqnarray*}
&&I_{1}\leq C\int_{\xi=\sum\limits_{j=1}^{k+1}\xi_{j}}\int_{\tau=\sum\limits_{j=1}^{k+1}\tau_{j}}
\frac{f(\prod\limits_{j=1}^{k+1}f_{j})|\phi^{\prime}(\xi)
-\phi^{\prime}(\xi_{k+1})|^{\frac{1}{2}-2\epsilon}|\xi_{1}|^{-\frac{1}{2}-4\epsilon}
|\xi_{2}||\xi_{3}|^{\frac{1}{4}+\epsilon}}
{\langle\sigma\rangle^{\frac{1}{2}-\frac{\epsilon}{12}}
\prod\limits_{j=1}^{k+1}\langle\sigma_{j}\rangle^{b}}d\delta\nonumber\\&&\leq C
\left\|I^{\frac{1}{2}-2\epsilon}(F_{k+1},F)\right\|_{L_{xt}^{2}}
\left\|D_{x}^{-\frac{1}{2}-4\epsilon}F_{1}\right\|_{L_{xt}^{\infty}}
\left(\prod_{j=4}^{k}\|F_{j}\|_{L_{xt}^{\infty}}\right)
\|D_{x}F_{2}\|_{L_{x}^{\infty}L_{t}^{2}}\|D_{x}^{\frac{1}{4}+\epsilon}F_{3}\|_{L_{x}^{2}L_{t}^{\infty}}\nonumber\\
&&\leq C\|F\|_{X_{0,\frac{1}{2}-\frac{\epsilon}{12}}}\prod\limits_{j=1}^{k+1}\|F_{j}\|_{X_{0,b}}\leq
C\left(\prod_{j=1}^{k+1}\|f_{j}\|_{L_{\xi\tau}^{2}}\right)\|f\|_{L_{\xi\tau}^{2}}.
\end{eqnarray*}
When (\ref{3.09})   is valid,  we  have  $|\xi|\sim|\xi_{k+1}|^{-1}$, thus,  we have
\begin{eqnarray}
&&\hspace{-1cm}K(\xi_{1},\xi_{2},\cdot\cdot\cdot,\xi_{k},\xi,\tau_{1},
\tau_{2},\cdot\cdot\cdot,\tau_{k},\tau)\leq C
\frac{|\xi_{1}|^{-\frac{1}{2}-4\epsilon}|\xi_{2}||\xi_{3}|^{\frac{1}{4}+\epsilon}
|\xi_{k+1}|^{\frac{1}{4}}|\xi|^{1-2\epsilon}}
{\langle\sigma\rangle^{\frac{1}{2}-\frac{\epsilon}{12}}
\prod\limits_{j=1}^{k+1}\langle\sigma_{j}\rangle^{b}}
.\label{3.017}
\end{eqnarray}
By using  (\ref{3.017}), the Plancherel identity,  the H\"older inequality, Lemma 2.11,
(\ref{2.08})-(\ref{2.010}), (\ref{2.056}),  (\ref{2.057}),
we have
\begin{eqnarray*}
&&I_{1}\leq C\int_{\xi=\sum\limits_{j=1}^{k+1}\xi_{j}}\int_{\tau=\sum\limits_{j=1}^{k+1}\tau_{j}}
\frac{f(\prod\limits_{j=1}^{k+1}f_{j})|\xi_{1}|^{-\frac{1}{2}-4\epsilon}
|\xi_{2}||\xi_{3}|^{\frac{1}{4}+\epsilon}
|\xi_{k+1}|^{\frac{1}{4}}|\xi|^{1-2\epsilon}}
{\langle\sigma\rangle^{\frac{1}{2}-\frac{\epsilon}{12}}
\prod\limits_{j=1}^{k+1}\langle\sigma_{j}\rangle^{b}}d\delta\nonumber\\
&&\leq C
\|D_{x}F_{2}\|_{L_{x}^{\infty}L_{t}^{2}}\|D_{x}^{\frac{1}{4}+\epsilon}F_{3}\|_{L_{x}^{2}L_{t}^{\infty}}
\|D_{x}^{-\frac{1}{2}-4\epsilon}F_{1}\|_{L_{xt}^{\infty}}
\left(\prod_{j=4}^{k}\|F_{j}\|_{L_{xt}^{\infty}}\right)\nonumber\\&&\qquad\qquad\times
\|D_{x}^{\frac{1}{4}}F_{k+1}\|_{L_{x}^{\frac{2}{1-2\epsilon}}L_{t}^{\infty}}
\|D_{x}^{1-2\epsilon}F\|_{L_{x}^{\frac{1}{\epsilon}}L_{t}^{2}}\nonumber\\
&&\leq C\|F\|_{X_{0,\frac{1}{2}-\frac{\epsilon}{12}}}\prod\limits_{j=1}^{k+1}\|F_{j}\|_{X_{0,b}}\leq
C\left(\prod_{j=1}^{k+1}\|f_{j}\|_{L_{\xi\tau}^{2}}\right)\|f\|_{L_{\xi\tau}^{2}}.
\end{eqnarray*}
(4)When $(\xi_{1},\xi_{2},\cdot\cdot\cdot,\xi_{k},\xi,\tau_{1},\tau_{2},
\cdot\cdot\cdot,\tau_{k},\tau)\in \Omega_{3}$, we consider
\begin{eqnarray}
&&|1+\frac{\gamma}{3\beta\xi_{3}^{2}\xi_{4}^{2}}|\geq\frac{1}{2},\label{3.018}\\
&&|1+\frac{\gamma}{3\beta\xi_{3}^{2}\xi_{4}^{2}}|<\frac{1}{2},\label{3.019}
\end{eqnarray}
respectively.

\noindent When  (\ref{3.018}) is valid, we have
\begin{eqnarray}
K(\xi_{1},\xi_{2},\cdot\cdot\cdot,\xi_{k},\xi,\tau_{1},\tau_{2},\cdot\cdot\cdot,\tau_{k},\tau)\leq C
\frac{|\xi_{3}^{2}-\xi_{4}^{2}|^{\frac{1}{2}}}
{\langle\sigma\rangle^{\frac{1}{2}-\frac{\epsilon}{12}}\prod\limits_{j=2}^{k+1}
\langle\xi_{j}\rangle^{s}\prod\limits_{j=1}^{k+1}\langle\sigma_{j}\rangle^{b}}\nonumber\\
\leq C
\frac{|\phi^{\prime}(\xi_{3})-\phi^{\prime}(\xi_{4})|^{\frac{1}{2}}}
{\langle\sigma\rangle^{\frac{1}{2}-\frac{\epsilon}{12}}\prod\limits_{j=2}^{k+1}
\langle\xi_{j}\rangle^{s}\prod\limits_{j=1}^{k+1}\langle\sigma_{j}\rangle^{b}}
.\label{3.020}
\end{eqnarray}
 By using (\ref{3.020}),  the Plancherel identity, the  H\"older inequality and
  (\ref{2.04}),  (\ref{2.06}),  (\ref{2.07}) as well as Lemma 2.5,
we have
\begin{eqnarray*}
&&I_{1}\leq C\int_{\xi=\sum\limits_{j=1}^{k+1}\xi_{j}}\int_{\tau=\sum\limits_{j=1}^{k+1}\tau_{j}}
\frac{f|\phi^{\prime}(\xi_{3})-\phi^{\prime}(\xi_{4})|^{\frac{1}{2}}
\prod\limits_{j=2}^{k+1}\langle\xi_{j}\rangle^{-s}
(\prod\limits_{j=1}^{k+1}f_{j})}{\langle\sigma\rangle^{\frac{1}{2}-\frac{\epsilon}{12}}
\prod\limits_{j=1}^{k+1}\langle\sigma_{j}\rangle^{b}}d\delta\nonumber\\&&\leq C
\left\|F_{1}\right\|_{L_{xt}^{8}}\left\|F_{2}\right\|_{L_{xt}^{8}}\|F\|_{L_{xt}^{\frac{8}{1+\epsilon}}}
\|I^{\frac{1}{2}}(F_{3},F_{4})\|_{L_{xt}^{2}}
\prod\limits_{j=5}^{k+1}\|J_{x}^{-\frac{k-4}{2(k-3)}-\epsilon}F_{j}\|_{L_{xt}^{\frac{8(k-3)}{1-\epsilon}}}\nonumber\\
&&\leq C\left(\prod_{j=1}^{k+1}\|F_{j}\|_{X_{0,b}}\right)
\|F\|_{X_{0,\frac{1}{2}-\frac{\epsilon}{12}}}\leq C\left(\prod_{j=1}^{k+1}
\|f_{j}\|_{L_{\xi\tau}^{2}}\right)\|f\|_{L_{\xi\tau}^{2}}.
\end{eqnarray*}
When (\ref{3.019})  is valid,  we have $|\xi_{3}|\sim|\xi_{4}|^{-1}$.

\noindent We consider $|\xi|\leq a,|\xi|\geq a,$ respectively.
\noindent When $|\xi|\leq a,$  we have
\begin{eqnarray}
&&K(\xi_{1},\xi_{2},\cdot\cdot\cdot,\xi_{k},\xi,\tau_{1},\tau_{2},\cdot\cdot\cdot,\tau_{k},\tau)\leq C
\frac{\prod\limits_{j=1}^{k} \langle \xi_{j}\rangle^{-s}}
{\langle\sigma\rangle^{\frac{1}{2}-\frac{\epsilon}{12}}\prod\limits_{j=1}^{k+1}\langle\sigma_{j}\rangle^{b}}
.\label{3.021}
\end{eqnarray}
By using  (\ref{3.021}), the Plancherel identity,  the H\"older inequality, (\ref{2.04}) and (\ref{2.056}),
we have
\begin{eqnarray*}
&&I_{1}\leq C\int_{\xi=\sum\limits_{j=1}^{k+1}\xi_{j}}\int_{\tau=\sum\limits_{j=1}^{k+1}\tau_{j}}
\frac{\prod\limits_{j=1}^{k} \langle \xi_{j}\rangle^{-s}f(\prod\limits_{j=1}^{k+1}f_{j})}
{\langle\sigma\rangle^{\frac{1}{2}-\frac{\epsilon}{12}}
\prod\limits_{j=1}^{k+1}\langle\sigma_{j}\rangle^{b}}d\delta\nonumber\\&&\leq C
\left(\prod_{j=1}^{k}\|J_{x}^{-s}F_{j}\|_{L_{xt}^{2k}}\right)\|F_{k+1}\|_{L_{xt}^{\infty}}
\|F\|_{L_{xt}^{2}}\nonumber\\
&&\leq C\|F\|_{X_{0,\frac{1}{2}-\frac{\epsilon}{12}}}\prod\limits_{j=1}^{k+1}\|F_{j}\|_{X_{0,b}}\leq
C\left(\prod_{j=1}^{k+1}\|f_{j}\|_{L_{\xi\tau}^{2}}\right)\|f\|_{L_{\xi\tau}^{2}}.
\end{eqnarray*}
When $|\xi|\geq a,$ we consider
\begin{eqnarray}
&&|1+\frac{\gamma}{3\beta\xi_{3}^{2}\xi_{5}^{2}}|\geq\frac{1}{2},\label{3.022}\\
&&|1+\frac{\gamma}{3\beta\xi_{3}^{2}\xi_{5}^{2}}|<\frac{1}{2},\label{3.023}
\end{eqnarray}
respectively.

\noindent When (\ref{3.022}) is valid,  this case can be proved similarly to (\ref{3.018}).
\noindent When (\ref{3.023}) is valid, we have that $|\xi_{3}|\sim |\xi_{4}|^{-1}\sim |\xi_{5}|^{-1},$ we consider
(\ref{3.08}),  (\ref{3.09}),  respectively.

\noindent When  (\ref{3.08}) is valid, we have
\begin{eqnarray}
K(\xi_{1},\xi_{2},\cdot\cdot\cdot,\xi_{k},\xi,\tau_{1},\tau_{2},\cdot\cdot\cdot,\tau_{k},\tau)\leq C
\frac{|\xi|^{1-4\epsilon}|\xi_{3}||\xi_{4}|^{-\frac{1}{4}+\epsilon}|\xi_{5}|^{\frac{1}{4}+\epsilon}
\prod\limits_{j=1}^{2}|\xi_{j}|^{-\frac{1}{2}-4\epsilon}}
{\langle\sigma\rangle^{\frac{1}{2}-\frac{\epsilon}{12}}\prod\limits_{j=1}^{k+1}\langle\sigma_{j}\rangle^{b}}.\label{3.024}
\end{eqnarray}
By using  (\ref{3.024}), the Plancherel identity,  the H\"older inequality,
Lemmas 2.5, 2.11, (\ref{2.08}) and (\ref{2.09}),  (\ref{2.056}),
we have
\begin{eqnarray*}
&&I_{1}\leq C\int_{\xi=\sum\limits_{j=1}^{k+1}\xi_{j}}\int_{\tau=\sum\limits_{j=1}^{k+1}\tau_{j}}
\frac{|\xi|^{1-4\epsilon}|\xi_{3}||\xi_{4}|^{\frac{1}{4}+\epsilon}|\xi_{5}|^{-\frac{1}{4}+\epsilon}
\prod\limits_{j=1}^{2}|\xi_{j}|^{-\frac{1}{2}-4\epsilon}f(\prod\limits_{j=1}^{k+1}f_{j})}
{\langle\sigma\rangle^{\frac{1}{2}-\frac{\epsilon}{12}}
\prod\limits_{j=1}^{k+1}\langle\sigma_{j}\rangle^{b}}d\delta\nonumber\\&&\leq C
\left[\prod\limits_{j=1}^{2}\left\|D_{x}^{-\frac{1}{2}-4\epsilon}F_{j}\right\|_{L_{xt}^{\infty}}\right]
\|D_{x}F_{3}\|_{L_{x}^{\infty}L_{t}^{2}}
\|D_{x}^{\frac{1}{4}+\epsilon}F_{4}\|_{L_{x}^{2}L_{t}^{\infty}}\|D_{x}^{-\frac{1}{4}+\epsilon}F_{5}\|_{L_{xt}^{\infty}}
\|I^{\frac{1}{2}-2\epsilon}(F,F_{k+1})\|_{L_{xt}^{2}}\nonumber\\&&\qquad\qquad\times
\prod\limits_{j=6}^{k}\|F_{j}\|_{L_{xt}^{\infty}}\nonumber\\
&&\leq C\|F\|_{X_{0,\frac{1}{2}-\frac{\epsilon}{12}}}\prod\limits_{j=1}^{k+1}\|F_{j}\|_{X_{0,b}}\leq
C\left(\prod_{j=1}^{k+1}\|f_{j}\|_{L_{\xi\tau}^{2}}\right)\|f\|_{L_{\xi\tau}^{2}}.
\end{eqnarray*}
\noindent When  (\ref{3.09}) is valid, we have
\begin{eqnarray}
K(\xi_{1},\xi_{2},\cdot\cdot\cdot,\xi_{k},\xi,\tau_{1},\tau_{2},\cdot\cdot\cdot,\tau_{k},\tau)\leq C
\frac{|\xi_{k+1}|^{\frac{1}{4}}|\xi|^{1-2\epsilon}
\prod\limits_{j=1}^{3}|\xi_{j}|^{\frac{1}{6}}}
{\langle\sigma\rangle^{\frac{1}{2}-\frac{\epsilon}{12}}\prod\limits_{j=1}^{k+1}\langle\sigma_{j}\rangle^{b}}.\label{3.025}
\end{eqnarray}
By using  (\ref{3.025}), the Plancherel identity,  the H\"older inequality,
 (\ref{2.05}) and (\ref{2.010}),  (\ref{2.057}),
we have
\begin{eqnarray*}
&&I_{1}\leq C\int_{\xi=\sum\limits_{j=1}^{k+1}\xi_{j}}\int_{\tau=\sum\limits_{j=1}^{k+1}\tau_{j}}
\frac{|\xi_{k+1}|^{\frac{1}{4}}|\xi|^{1-2\epsilon}
\prod\limits_{j=1}^{3}|\xi_{j}|^{\frac{1}{6}}f(\prod\limits_{j=1}^{k+1}f_{j})}
{\langle\sigma\rangle^{\frac{1}{2}-\frac{\epsilon}{12}}
\prod\limits_{j=1}^{k+1}\langle\sigma_{j}\rangle^{b}}d\delta\nonumber\\&&\leq C
\left[\prod\limits_{j=1}^{3}\left\|D_{x}^{\frac{1}{6}}F_{j}\right\|_{L_{xt}^{6}}\right]
\left[\prod\limits_{j=4}^{k}\|F_{j}\|_{L_{xt}^{\infty}}\right]
\|D_{x}^{\frac{1}{4}}F_{k+1}\|_{L_{x}^{\frac{2}{1-2\epsilon}}L_{t}^{\infty}}
\|D_{x}^{1-2\epsilon}F\|_{L_{x}^{\frac{1}{\epsilon}}L_{t}^{2}}
\nonumber\\
&&\leq C\|F\|_{X_{0,\frac{1}{2}-\frac{\epsilon}{12}}}\prod\limits_{j=1}^{k+1}\|F_{j}\|_{X_{0,b}}\leq
C\left(\prod_{j=1}^{k+1}\|f_{j}\|_{L_{\xi\tau}^{2}}\right)\|f\|_{L_{\xi\tau}^{2}}.
\end{eqnarray*}
(5)Case$(\xi_{1},\xi_{2},\cdot\cdot\cdot,\xi_{k},\xi,\tau_{1},\tau_{2},
\cdot\cdot\cdot,\tau_{k},\tau)\in \Omega_{4}$,
this case can be proved similarly to Case (4).\\
\noindent
(6)Case $(\xi_{1},\xi_{2},\cdot\cdot\cdot,\xi_{k},\xi,\tau_{1},\tau_{2},
\cdot\cdot\cdot,\tau_{k},\tau)\in \Omega_{5}$,
we consider
\begin{eqnarray}
&&|1+\frac{\gamma}{3\beta\xi_{5}^{2}\xi_{6}^{2}}|\geq\frac{1}{2},\label{3.026}\\
&&|1+\frac{\gamma}{3\beta\xi_{5}^{2}\xi_{6}^{2}}|<\frac{1}{2},\label{3.027}
\end{eqnarray}
respectively.

\noindent When (\ref{3.026}) is valid, we consider $k=5$, $k\geq6,$  respectively.

\noindent When $k=5,$  we have
\begin{eqnarray}
K(\xi_{1},\xi_{2},\cdot\cdot\cdot,\xi_{k},\xi,\tau_{1},\tau_{2},\cdot\cdot\cdot,\tau_{k},\tau)\leq C
\frac{\langle\xi_{6}\rangle^{-\frac{1}{3}-3\epsilon}\prod\limits_{j=1}^{5}|\xi_{j}|^{\frac{1}{6}}}
{\langle\sigma\rangle^{\frac{1}{2}-\frac{\epsilon}{12}}
\prod\limits_{j=1}^{6}\langle\sigma_{j}\rangle^{b}}
.\label{3.028}
\end{eqnarray}
By using (\ref{3.028}), the Plancherel identity, the  H\"older inequality
and   (\ref{2.04}), (\ref{2.05}) and (\ref{2.07}),
we have
\begin{eqnarray*}
&&I_{1}\leq C\int_{\xi=\sum\limits_{j=1}^{6}\xi_{j}}\int_{\tau=\sum\limits_{j=1}^{6}\tau_{j}}
\frac{ff_{6}(\prod\limits_{j=1}^{5}|\xi_{j}|^{\frac{1}{6}}f_{j})\langle\xi_{6}\rangle^{-\frac{1}{3}-3\epsilon}}
{\langle\sigma\rangle^{\frac{1}{2}-\frac{\epsilon}{12}}
\prod\limits_{j=1}^{6}\langle\sigma_{j}\rangle^{b}}d\delta\nonumber\\&&\leq C
\|F\|_{L_{xt}^{\frac{8}{1+\epsilon}}}\|J_{x}^{-\frac{1}{3}-3\epsilon}F_{6}\|_{L_{xt}^{\frac{24}{1-3\epsilon}}}
\prod\limits_{j=1}^{5}\left\|D_{x}^{\frac{1}{6}}F_{j}\right\|_{L_{xt}^{6}}
\nonumber\\&&\leq
C\left(\prod_{j=1}^{6}\|f_{j}\|_{L_{xt}^{2}}\right)\|f\|_{L_{xt}^{2}}
.
\end{eqnarray*}
When $k\geq6,$
we have
\begin{eqnarray}
K(\xi_{1},\xi_{2},\cdot\cdot\cdot,\xi_{k},\xi,\tau_{1},\tau_{2},\cdot\cdot\cdot,\tau_{k},\tau)\leq C
\frac{|\xi_{5}^{2}-\xi_{6}^{2}|^{\frac{1}{2}}\prod\limits_{j=1}^{4}\langle\xi_{j}
\rangle^{-s}\prod\limits_{j=7}^{k+1}\langle\xi_{j}\rangle^{-\frac{(k-4)s}{(k-5)}}}
{\langle\sigma\rangle^{\frac{1}{2}-\frac{\epsilon}{12}}
\prod\limits_{j=1}^{k+1}\langle\sigma_{j}\rangle^{b}}\nonumber\\
\leq C
\frac{|\phi^{\prime}(\xi_{5})-\phi^{\prime}(\xi_{6})|^{\frac{1}{2}}\prod\limits_{j=1}^{4}\langle\xi_{j}
\rangle^{-s}\prod\limits_{j=7}^{k+1}\langle\xi_{j}\rangle^{-\frac{(k-4)s}{(k-5)}}}
{\langle\sigma\rangle^{\frac{1}{2}-\frac{\epsilon}{12}}\prod\limits_{j=1}^{k+1}\langle\sigma_{j}\rangle^{b}}
.\label{3.028}
\end{eqnarray}
By using (\ref{3.028}),   the Plancherel identity, the  H\"older inequality
and   (\ref{2.04}),  (\ref{2.07}), Lemma 2.5,
we have
\begin{eqnarray*}
&&I_{1}\leq C\int_{\xi=\sum\limits_{j=1}^{k+1}\xi_{j}}\int_{\tau=\sum\limits_{j=1}^{k+1}\tau_{j}}
\frac{|\phi^{\prime}(\xi_{5})-\phi^{\prime}(\xi_{6})|^{\frac{1}{2}}
f(\prod\limits_{j=1}^{k+1}f_{j})\prod\limits_{j=1}^{4}\langle\xi_{j}
\rangle^{-s}\prod\limits_{j=7}^{k+1}\langle\xi_{j}\rangle^{-\frac{(k-4)s}{(k-5)}}}
{\langle\sigma\rangle^{\frac{1}{2}-\frac{\epsilon}{12}}\prod\limits_{j=1}^{k+1}
\langle\sigma_{j}\rangle^{b}}d\delta\nonumber\\&&\leq C
\|F\|_{L_{xt}^{\frac{8}{1+\epsilon}}}
\left(\prod\limits_{j=1}^{4}\left\|J_{x}^{-s}F_{j}\right\|_{L_{xt}^{2k}}\right)
\|I^{\frac{1}{2}}(F_{5},F_{6})\|_{L_{xt}^{2}}
\prod\limits_{j=7}^{k+1}\|J_{x}^{-\frac{(k-4)s}{(k-5)}}F_{j}\|_{L_{xt}^{\frac{8k(k-5)}{3k-16-k\epsilon}}}
\nonumber\\
&&\leq C\left(\prod_{j=1}^{k+1}\|F_{j}\|_{X_{0,b}}\right)
\|F\|_{X_{0,\frac{1}{2}-\frac{\epsilon}{12}}}\leq C
\left(\prod_{j=1}^{k+1}\|f_{j}\|_{L_{xt}^{2}}\right)\|f\|_{L_{xt}^{2}}.
\end{eqnarray*}
When  (\ref{3.027}) is valid, we have $|\xi_{5}|\sim|\xi_{6}|^{-1}$.
 If there exists some $k\in N,7\leq k\leq k+1$  such  that
 \begin{eqnarray*}
&&|1+\frac{\gamma}{3\beta\xi_{5}^{2}\xi_{k}^{2}}|\geq\frac{1}{2},
\end{eqnarray*}
 we can
use a proof similar to (\ref{3.026}) to derive the result, otherwise
 we have that $|\xi_{5}|\sim |\xi_{6}|^{-1}\sim |\xi_{k+1}|^{-1}$.

\noindent In this case,  we have
\begin{eqnarray}
&&K(\xi_{1},\xi_{2},\cdot\cdot\cdot,\xi_{k},\xi,\tau_{1},\tau_{2},\cdot\cdot\cdot,\tau_{k},\tau)\leq C
\frac{|\xi|^{1-2\epsilon}|\xi_{1}|\prod\limits_{j=2}^{5}|\xi_{j}|^{-\frac{1}{4}-2\epsilon}}
{\langle\sigma\rangle^{\frac{1}{2}-\frac{\epsilon}{12}}\prod\limits_{j=1}^{k+1}\langle\sigma_{j}\rangle^{b}}
.\label{3.029}
\end{eqnarray}
By using  (\ref{3.029}), the Plancherel identity,  the H\"older inequality,
 (\ref{2.08}),  (\ref{2.010}),   (\ref{2.045})  and   (\ref{2.056}),
we have
\begin{eqnarray*}
&&I_{1}\leq C\int_{\xi=\sum\limits_{j=1}^{k+1}\xi_{j}}
\int_{\tau=\sum\limits_{j=1}^{k+1}\tau_{j}}
\frac{f(\prod\limits_{j=1}^{k+1}f_{j})|\xi|^{1-2\epsilon}|\xi_{1}|
\prod\limits_{j=2}^{5}|\xi_{j}|^{-\frac{1}{4}-2\epsilon}}
{\langle\sigma\rangle^{\frac{1}{2}-\frac{\epsilon}{12}}
\prod\limits_{j=1}^{k+1}\langle\sigma_{j}\rangle^{b}}d\delta\nonumber\\
&&\leq C\|D_{x}F_{1}\|_{L_{x}^{\infty}L_{t}^{2}}
\left(\prod_{j=2}^{5}\|D_{x}^{-\frac{1}{4}-2\epsilon}F_{j}
\|_{L_{x}^{\frac{4}{1-\epsilon}}L_{t}^{\infty}}\right)
\left(\prod_{j=6}^{k+1}\|F_{j}\|_{L_{xt}^{\infty}}\right)
\|D_{x}^{1-2\epsilon}F\|_{L_{x}^{\frac{1}{\epsilon}}L_{t}^{2}}\nonumber\\
&&\leq C\|F\|_{X_{0,\frac{1}{2}-\frac{\epsilon}{12}}}
\prod\limits_{j=1}^{k+1}\|F_{j}\|_{X_{0,b}}\leq
C\left(\prod_{j=1}^{k+1}\|f_{j}\|_{L_{xt}^{2}}\right)
\|f\|_{L_{xt}^{2}}.
\end{eqnarray*}
(7)Case $(\xi_{1},\xi_{2},\cdot\cdot\cdot,\xi_{k},\xi,\tau_{1},
\tau_{2},\cdot\cdot\cdot,\tau_{k},\tau)\in \Omega_{6}$,
this case can be proved similarly to
$(\xi_{1},\xi_{2},\cdot\cdot\cdot,\xi_{k},\xi,\tau_{1},\tau_{2},
\cdot\cdot\cdot,\tau_{k},\tau)\in \Omega_{5}$.

\noindent(8)Case $(\xi_{1},\xi_{2},\cdot\cdot\cdot,\xi_{k},\xi,\tau_{1},
\tau_{2},\cdot\cdot\cdot,\tau_{k},\tau)\in \Omega_{7}$, we consider
$|\xi|\leq a,|\xi|\geq a$, respectively.

\noindent When $|\xi|\leq a,$  we consider (\ref{3.08}),  (\ref{3.09}), respectively.

\noindent When (\ref{3.08})  is valid,  we have
\begin{eqnarray}
&&K(\xi_{1},\xi_{2},\cdot\cdot\cdot,\xi_{k},\xi,\tau_{1},\tau_{2},\cdot\cdot\cdot,\tau_{k},\tau)\leq C
\frac{|\xi_{k+1}|^{1-4\epsilon}\prod\limits_{j=1}^{k+1}\langle\xi_{j}\rangle^{-s}}
{\langle\sigma\rangle^{\frac{1}{2}-\frac{\epsilon}{12}}\prod\limits_{j=1}^{k+1}\langle\sigma_{j}\rangle^{b}}
.\label{3.030}
\end{eqnarray}
By using (\ref{3.030}), the Plancherel identity, the  H\"older inequality
and   (\ref{2.04}), Lemma 2.5,
we have
\begin{eqnarray*}
&&I_{1}\leq C\int_{\xi=\sum\limits_{j=1}^{k+1}\xi_{j}}\int_{\tau=\sum\limits_{j=1}^{k+1}\tau_{j}}
\frac{|\xi_{k+1}|^{1-4\epsilon}f(\prod\limits_{j=1}^{k+1}f_{j})\prod\limits_{j=1}^{k+1}\langle\xi_{j}\rangle^{-s}}
{\langle\sigma\rangle^{\frac{1}{2}-\frac{\epsilon}{12}}
\prod\limits_{j=1}^{k+1}\langle\sigma_{j}\rangle^{b}}d\delta\nonumber\\&&\leq C
\left[\prod\limits_{j=1}^{k}\|D_{x}^{-s}F_{j}\|_{L_{xt}^{2k}}\right]
\|I^{\frac{1}{2}-2\epsilon}(F_{k+1},F)\|_{L_{xt}^{2}}
\nonumber\\
&&\leq C\left(\prod_{j=1}^{k+1}\|F_{j}\|_{X_{0,b}}\right)
\|F\|_{X_{0,\frac{1}{2}-\frac{\epsilon}{12}}}\leq C\left(\prod_{j=1}^{k+1}\|f_{j}\|_{L_{\xi\tau}^{2}}\right)\|f\|_{L_{\xi\tau}^{2}}.
\end{eqnarray*}
When (\ref{3.09}) is valid,  we consider $|\sigma|\geq |\xi_{1}|, |\sigma|\leq |\xi_{1}|,$ respectively.

\noindent When $|\sigma|\geq |\xi_{1}|,$
 we have
\begin{eqnarray}
&&K(\xi_{1},\xi_{2},\cdot\cdot\cdot,\xi_{k},\xi,\tau_{1},\tau_{2},\cdot\cdot\cdot,\tau_{k},\tau)\leq C
\frac{|\xi_{1}|^{-\frac{1}{2}+4\epsilon}\prod\limits_{j=1}^{k+1}\langle\xi_{j}\rangle^{-s}}
{\prod\limits_{j=1}^{k+1}\langle\sigma_{j}\rangle^{b}}
.\label{3.031}
\end{eqnarray}
By using (\ref{3.031}), the Plancherel identity, the  H\"older inequality
and   (\ref{2.04}),
we have
\begin{eqnarray*}
&&I_{1}\leq C\int_{\xi=\sum\limits_{j=1}^{k+1}\xi_{j}}\int_{\tau=\sum\limits_{j=1}^{k+1}\tau_{j}}
\frac{f(\prod\limits_{j=1}^{k+1}f_{j})\prod\limits_{j=1}^{k+1}\langle\xi_{j}\rangle^{-s-\frac{\frac{1}{2}-4\epsilon}{k+1}}}
{\prod\limits_{j=1}^{k+1}\langle\sigma_{j}\rangle^{b}}d\delta\nonumber\\&&\leq C
\left[\prod\limits_{j=1}^{k+1}\|D_{x}^{-s-\frac{\frac{1}{2}-4\epsilon}{k+1}}F_{j}\|_{L_{xt}^{2(k+1)}}\right]
\|\mathscr{F}^{-1}f\|_{L_{xt}^{2}}
\nonumber\\
&&\leq C\left(\prod_{j=1}^{k+1}\|F_{j}\|_{X_{0,b}}\right)
\|f\|_{L_{\xi\tau}^{2}}\leq C\left(\prod_{j=1}^{k+1}\|f_{j}\|_{L_{\xi\tau}^{2}}\right)\|f\|_{L_{\xi\tau}^{2}}.
\end{eqnarray*}
\noindent When $|\sigma|\leq |\xi_{1}|,$
 we have
\begin{eqnarray}
&&K(\xi_{1},\xi_{2},\cdot\cdot\cdot,\xi_{k},\xi,\tau_{1},\tau_{2},\cdot\cdot\cdot,\tau_{k},\tau)\leq C
\frac{|\xi_{1}|^{\frac{1}{2}}\prod\limits_{j=1}^{k+1}\langle\xi_{j}\rangle^{-s}}
{\langle\sigma\rangle^{\frac{1}{2}-\frac{\epsilon}{12}}\prod\limits_{j=1}^{k+1}\langle\sigma_{j}\rangle^{b}}
.\label{3.032}
\end{eqnarray}
By using (\ref{3.032}), the Plancherel identity, the  H\"older inequality
and   (\ref{2.04}), (\ref{2.010}), (\ref{2.057}),
we have
\begin{eqnarray*}
&&I_{1}\leq C\int_{\xi=\sum\limits_{j=1}^{k+1}\xi_{j}}\int_{\tau=\sum\limits_{j=1}^{k+1}\tau_{j}}
\frac{|\xi_{1}|^{\frac{1}{2}}f(\prod\limits_{j=1}^{k+1}f_{j})\prod\limits_{j=1}^{k+1}\langle\xi_{j}\rangle^{-s}}
{\langle\sigma\rangle^{\frac{1}{2}-\frac{\epsilon}{12}}
\prod\limits_{j=1}^{k+1}\langle\sigma_{j}\rangle^{b}}d\delta\nonumber\\&&\leq C
\left[\prod\limits_{j=1}^{k}\|D_{x}^{-s}F_{j}\|_{L_{xt}^{2k}}\right]\|D_{x}^{1-2\epsilon}F_{k+1}\|_{L_{x}^{\frac{1}{\epsilon}}L_{t}^{2}}
\|D_{x}^{\frac{1}{4}}F\|_{L_{x}^{\frac{2}{1-2\epsilon}}L_{t}^{\infty}}
\nonumber\\
&&\leq C\left(\prod_{j=1}^{k+1}\|F_{j}\|_{X_{0,b}}\right)
\|F\|_{X_{0,\frac{1}{2}-\frac{\epsilon}{12}}}\leq C\left(\prod_{j=1}^{k+1}\|f_{j}\|_{L_{xt}^{2}}\right)\|f\|_{L_{xt}^{2}}.
\end{eqnarray*}
When $|\xi|\geq a,$  we have
\begin{eqnarray}
&&K(\xi_{1},\xi_{2},\cdot\cdot\cdot,\xi_{k},\xi,\tau_{1},\tau_{2},\cdot\cdot\cdot,\tau_{k},\tau)\leq C
\frac{|\xi|^{\frac{3-\epsilon}{18}}|\xi_{1}|^{\frac{1}{6}}|\xi_{2}|^{\frac{1}{6}}
\prod\limits_{j=3}^{k+1}\langle\xi_{j}\rangle^{\frac{1-2ks}{2(k-1)}}}
{\langle\sigma\rangle^{\frac{1}{2}-\frac{\epsilon}{12}}\prod\limits_{j=1}^{k+1}\langle\sigma_{j}\rangle^{b}}
.\label{3.033}
\end{eqnarray}
By using (\ref{3.033}), the Plancherel identity, the  H\"older inequality
and   (\ref{2.04}), (\ref{2.05}),  (\ref{2.0888888}),
we have
\begin{eqnarray*}
&&I_{1}\leq C\int_{\xi=\sum\limits_{j=1}^{k+1}\xi_{j}}\int_{\tau=\sum\limits_{j=1}^{k+1}\tau_{j}}
\frac{|\xi|^{\frac{3-\epsilon}{18}}|\xi_{1}|^{\frac{1}{6}}|\xi_{2}|^{\frac{1}{6}}
\prod\limits_{j=3}^{k+1}\langle\xi_{j}\rangle^{\frac{1-2ks}{2(k-1)}}
f(\prod\limits_{j=1}^{k+1}f_{j})}
{\langle\sigma\rangle^{\frac{1}{2}-\frac{\epsilon}{12}}
\prod\limits_{j=1}^{k+1}\langle\sigma_{j}\rangle^{b}}d\delta\nonumber\\&&\leq C
\|D_{x}^{\frac{3-\epsilon}{18}}F\|_{L_{xt}^{\frac{6}{1+\epsilon}}}
\left[\prod\limits_{j=1}^{2}\|D_{x}^{\frac{1}{6}}F_{j}\|_{L_{xt}^{6}}\right]
\left[\prod\limits_{j=3}^{k+1}\|D_{x}^{\frac{1-2ks}{2(k-1)}}
F_{j}\|_{L_{xt}^{\frac{6(k-1)}{3-\epsilon}}}\right]
\nonumber\\
&&\leq C\left(\prod_{j=1}^{k+1}\|F_{j}\|_{X_{0,b}}\right)
\|F\|_{X_{0,\frac{1}{2}-\frac{\epsilon}{12}}}\leq C
\left(\prod_{j=1}^{k+1}\|f_{j}\|_{L_{\xi\tau}^{2}}\right)\|f\|_{L_{\xi\tau}^{2}}.
\end{eqnarray*}

We have completed the proof of Lemma 3.1.
\begin{Lemma}\label{Lemma3.2}
Let $s\geq\frac{1}{2}-\frac{2}{k}+2\epsilon,k\geq5$,   $b=\frac{1}{2}+\frac{\epsilon}{24}$
 and $b^{\prime}=-\frac{1}{2}+\frac{\epsilon}{12}$ and  $g=\psi(t)u$.
Then, we have
\begin{eqnarray}
\left\|g^{k+1}\right\|_{X_{s,b^{\prime}}}\leq C\|f\|_{X_{s,b}}^{k+1}.\label{3.034}
\end{eqnarray}
\end{Lemma}

Lemma 3.2 can be proved similarly to Lemma 3.1.

\begin{Lemma}\label{Lemma3.3}
Let $s\geq\frac{1}{2}-\frac{2}{k}+2\epsilon,k\geq5$,  $b=\frac{1}{2}+\frac{\epsilon}{24}$
 and $b^{\prime}=-\frac{1}{2}+\frac{\epsilon}{12}$ and   $g=\psi(t)u$.
Then, we have
\begin{eqnarray}
\left\|\partial_{x} (g^{k+1})\right\|_{\tilde{X}_{s,b^{\prime}}}\leq
 C\|g\|_{\tilde{X}_{s,b}}^{k+1}.\label{3.035}
\end{eqnarray}
\end{Lemma}
\noindent {\bf Proof.}Since $\|g\|_{\tilde{X}_{s,b}}=
\|g\|_{X_{s,b}}+\|\partial_{x}^{-1}g\|_{X_{s,b}}$,we have
\begin{eqnarray}
&&\|\partial_{x}(g^{k+1})\|_{\tilde{X}_{s,b^{\prime}}}=
\|\partial_{x}(g^{k+1})\|_{X_{s,b^{\prime}}}
+\|g^{k+1}\|_{X_{s,b^{\prime}}},\label{3.036}
\end{eqnarray}
using Lemma 3.1 and Lemma 3.2,  we have
\begin{eqnarray*}
&&\|\partial_{x}(g^{k+1})\|_{\tilde{X}_{s,b^{\prime}}}
\leq C\|g\|_{X_{s,b}}^{k+1}
\leq C\|g\|_{\tilde{X}_{s,b}}^{k+1}.
\end{eqnarray*}

\noindent We have completed the proof of Lemma 3.3.
\bigskip

\noindent {\large\bf 4. Proof of Theorem  1.1}

\setcounter{equation}{0}

 \setcounter{Theorem}{0}

\setcounter{Lemma}{0}

\setcounter{section}{4}
\noindent {\bf Proof.}
Obviously, (\ref{1.01})-(\ref{1.02}) are equivalent to
\begin{eqnarray*}
&&u(t)=U^{\gamma,\beta}(t)u_{0}-\frac{1}{k+1}\int_{0}^{t}
U^{\gamma,\beta}(t-t^{\prime})\partial_{x}[(\psi(\tau)u)^{k+1}]d\tau.
\end{eqnarray*}
For $u_{0}\in H^{s}(\R)$ and $\delta\in(0,1]$, we define
\begin{eqnarray}
&&\Gamma (v)=\psi(t)U^{\gamma,\beta}(t)u_{0}+\frac{1}{k+1}\psi_{\delta}(t)\int_{0}^{t}
U^{\gamma,\beta}(t-\tau)\partial_{x}((\psi(\tau)v)^{k+1})d\tau.\label{4.01}
\end{eqnarray}
We define
$B(0,r)=\{u\in X_{s,b}\cap C([-\delta,\delta],H^{s}(\R)),\|u\|_{X_{s,b}}\leq r:=2C\|u_{0}\|_{H^{s}(\SR)}\}$.
By using  Lemmas  2.1, 3.1 and  choosing  sufficiently  small  $\delta>0$  such
 that
 \begin{eqnarray*}
C\delta^{b^{\prime}+1-b}(2C\|u_{0}\|_{H^{s}(\SR)})^{k+1}\leq C\|u_{0}\|_{H^{s}},
 \end{eqnarray*}
  we have
\begin{eqnarray}
&&\|\Gamma(v)\|_{X_{s,b}}\leq C\|u_{0}\|_{H^{s}}+C\delta^{b^{\prime}+1-b}
\|\partial_{x}(v^{k+1})\|_{X_{s,b^{\prime}}}\nonumber\\
&&\leq C\|u_{0}\|_{H^{s}}+C\delta^{b^{\prime}+1-b}\|v\|_{X_{s,b}}^{k+1}\nonumber\\
&&\leq C\|u_{0}\|_{H^{s}}+C\delta^{b^{\prime}+1-b}(2C\|u_{0}\|_{H^{s}(\SR)})^{k+1}
\leq 2C\|u_{0}\|_{H^{s}}.\label{4.02}
\end{eqnarray}
By a similar calculation, we have
\begin{eqnarray}
&&\|\Gamma(v_{1})-\Gamma(v_{2})\|_{X_{s,b}}\leq C\delta^{b^{\prime}+1-b}\|v_{1}-v_{2}\|_{X_{s,b}}(\|u\|_{X_{s,b}}^{k}+\|v\|_{X_{s,b}}^{k})\nonumber\\&&\leq \frac{1}{2}\|v_{1}-v_{2}\|_{X_{s,b}}.\label{4.03}
\end{eqnarray}
Thus $\Gamma$ is a contraction mapping from the closed ball
$$B(0,r)=\{u\in X_{s,b}\cap C([-\delta,\delta],H^{s}(\R)),\|u\|_{X_{s,b}}\leq r\}$$
into itself. From  the fixed point theorem and  (\ref{4.03}), we have $\Gamma(v)=v$.
The uniqueness of  solution to (\ref{4.01}) is esasily derived from (\ref{4.03}).

The rest of the local well-posedness results of Theorem 1.1 follow from a
 standard argument, for instance, see \cite{KPV1993}.

This completes the proof of Theorem 1.1.

\bigskip

\noindent {\large\bf 5. Proof of Theorem  1.2}

\setcounter{equation}{0}

 \setcounter{Theorem}{0}

\setcounter{Lemma}{0}

\setcounter{section}{5}

\bigskip

By using Lemmas 2.1, 3.3 and the fixed point argument as  well as a
  proof similar to Theorem 1.1, we obtain Theorem 1.2.

\bigskip

\noindent {\large\bf 6. Proof of Theorem  1.3}

\setcounter{equation}{0}

 \setcounter{Theorem}{0}

\setcounter{Lemma}{0}

\setcounter{section}{6}
\noindent {\bf Proof.}In this section, inspired \cite{BS,LV,LS,T},  by  we study the relationship between the
 solution to (\ref{1.01})-(\ref{1.02}) and the solution to
\begin{eqnarray}
     && v_{t}-\beta\partial_{x}^{3}v+\frac{1}{k+1}(v^{k+1})_{x}=0,\label{6.01}\\
     &&v(x,0)=v_{0}(x),\label{6.02}
 \end{eqnarray}
as $\gamma \rightarrow 0$.

\noindent From (\ref{1.01}), we have
 \begin{eqnarray}
      J_{x}^{s}u_{t}-\beta\partial_{x}^{3}J_{x}^{s}u-\gamma\partial_{x}^{-1}J_{x}^{s}u+\frac{1}{k+1}J_{x}^{s}(u^{k+1})_{x}=0,k\geq5\label{6.099}
    \end{eqnarray}
Multiplying    by  $J_{x}^{s}u$  on both sides of (\ref{6.099}) and integration by parts with respect to $x$ on $\R$
  as well as $H^{s-1}(\R)\hookrightarrow L^{\infty}$ with $s>\frac{3}{2},$ by using Lemma 2.12,
we obtain
\begin{eqnarray}
&&\frac{1}{2}\frac{d}{dt}\|u\|_{H^{s}}^{2}=\int_{\SR}J_{x}^{s}uJ_{x}^{s}u_{t}dx=\beta\int_{\SR}J_{x}^{s}u \partial_{x}^{3}J_{x}^{s}udx
+\gamma\int_{\SR}J_{x}^{s}u \partial_{x}^{-1}J_{x}^{s}udx\nonumber\\
&&-\int_{\SR}J_{x}^{s}u J_{x}^{s}(u^{k}u_{x})dx= -\int_{\SR}J_{x}^{s}u J_{x}^{s}(u^{k}u_{x})dx\nonumber\\
&&=-\int_{\SR}J_{x}^{s}u[J_{x}^{s},u^{k}]u_{x}dx-\int_{\SR}(J_{x}^{s}u)(u^{k}J_{x}^{s}u_{x})dx\nonumber\\
&&=-\int_{\SR}J_{x}^{s}u[J_{x}^{s},u^{k}]u_{x}dx+\frac{k}{2}\int_{\SR}u^{k-1}u_{x}(J_{x}^{s}u_{x})^{2}dx\nonumber\\
&&\leq C\left\|J_{x}^{s}u\right\|_{L^{2}}\left[\|u^{k-1}u_{x}\|_{L^{\infty}}
\|J_{x}^{s}u\|_{L^{2}}+\|J_{x}^{s}(u)^{k}\|_{L^{2}}\|u_{x}\|_{L^{\infty}}\right]\nonumber\\
&&+C\left\|u^{k-1}u_{x}\right\|_{L^{\infty}}\left\|J_{x}^{s}u\right\|_{L^{2}}^{2}\leq C_{0}\left\|u\right\|_{H^{s}}^{k+2}\leq C_{0}\left\|u\right\|_{X_{s}}^{k+2}.\label{6.03}
\end{eqnarray}
Here $C_{0}$ is a constant independent of $\gamma$. Similarly, we obtain
\begin{eqnarray}
&&\frac{1}{2}\frac{d}{dt}\|\partial_{x}^{-1}u\|_{H^{s}}^{2}=-2\int_{\SR}J_{x}^{s}\partial_{x}^{-1}uJ_{x}^{s}(u^{k+1})dx\nonumber\\
&&\leq C_{0}\|J_{x}^{s}\partial_{x}^{-1}u\|_{L^{2}}\|u^{k+1}\|_{H^{s}}\leq C_{0}\|u\|_{X_{s}}^{k+2}.\label{6.04}
\end{eqnarray}
Then, using (\ref{6.03}) and (\ref{6.04}),  we have
\begin{eqnarray}
&&\frac{d}{dt}\|u\|_{X_{s}}^{2}\leq C_{0}\|u\|_{X_{s}}^{k+2},\label{6.05}
\end{eqnarray}
from (\ref{6.05}), we have
\begin{eqnarray}
&&\frac{d}{dt}\|u\|_{X_{s}}\leq C_{0}\|u\|_{X_{s}}^{k+1}.\label{6.06}
\end{eqnarray}
When $t<\min\left\{T,\frac{1}{Ck\|u_{0}\|_{X_{s}}^{k}}\right\}$, where $T$  is the time lifespan of  the solution
to (\ref{1.01})-(\ref{1.02}) for data in $X_{s}(\R)$ with $s>\frac{3}{2}$ in Theorem 1.2,  by    using (\ref{6.06}), we have
\begin{eqnarray}
&&\|u\|_{X_{s}}\leq \frac{k^{\frac{1}{k}}\|u_{0}\|_{X_{s}}}{\sqrt{1-Ck\|u_{0}\|_{X_{s}}^{k}t}}.\label{6.07}
\end{eqnarray}
Let $u:=u^{\gamma}$ and  the solution to (\ref{1.01}).
 Therefore $w:=u-v$ satisfies the equation
\begin{eqnarray}
&w_{t}-\beta\partial_{x}^{3}w+\gamma\partial_{x}^{-1}u+\frac{1}{k+1}(w\sum\limits_{j=0}^{k}(v+w)^{j}v^{k-j})_{x}=0,k\geq5,\label{6.08}\\
         &w(x,0)=u_{0}(x)-v_{0}(x).\label{6.09}
\end{eqnarray}
Multiplying by  $w$ on  both sides of (\ref{6.08}) and  integrating by parts  with respect to $x$ on  $\R$, we obtain
\begin{eqnarray}
&&\frac{1}{2}\frac{d}{dt}\|w\|_{L^{2}}^{2}=-\frac{1}{k+1}\int_{\SR}w(w\sum\limits_{j=0}^{k}(v+w)^{j}v^{k-j})_{x}dx
+\int_{\SR}w\gamma\partial_{x}^{-1}u)dx\nonumber\\
&&=-\frac{1}{k+1}\int_{\SR}w(w\sum\limits_{j=0}^{k}u^{j}v^{k-j})_{x}dx+
\gamma\int_{\SR}w\partial_{x}^{-1}udx\nonumber\\
&&\leq C\sum\limits_{j=0}^{k}\|(u^{j}v^{k-j})_{x}\|_{L^{\infty}}\|w\|_{L^{2}}^{2}
+|\gamma|\|w\|_{L^{2}}\|\partial_{x}^{-1}u\|_{L^{2}}\nonumber\\
&&\leq C\left[\|u\|_{H^{s}}^{k}+\|v\|_{H^{s}}^{k}\right]\|w\|_{L^{2}}^{2}
+C|\gamma|\|w\|_{L^{2}}\|u\|_{X_{s}}.\label{6.010}
\end{eqnarray}
From (\ref{6.09}),  we have
\begin{eqnarray}
&&\frac{d}{dt}\|w\|_{L^{2}}\leq C\sup\limits_{t\in [0,T]}\left[\|u\|_{X_{s}}+\|v\|_{X_{s}}\right]^{k}\|w\|_{L^{2}}
+C|\gamma|\sup\limits_{t\in [0,T]}\|u\|_{X_{s}}.\label{6.011}
\end{eqnarray}
By using the Gronwall's inequality and (\ref{6.011}),  we can get
\begin{eqnarray}
\|w\|_{L^{2}}\leq e^{CT\sup\limits_{t\in [0,T]}\left[\|u\|_{X_{s}}+\|v\|_{X_{s}}\right]^{k}}
\left[\|u_{0}-v_{0}\|_{L^{2}}+C|\gamma| T\sup\limits_{t\in [0,T]}
\|u\|_{X_{s}}\right].\label{6.012}
\end{eqnarray}
Thus, when $\gamma\rightarrow 0$ and  $\|u_{0}-v_{0}\|_{L^{2}}\rightarrow0$, then, we have
$\|w\|_{L^{2}}=\|u-v\|_{L^{2}}\rightarrow0.$

 This completes the proof of Theorem 1.3.

\bigskip
\leftline{\large \bf Acknowledgments}

\noindent This work is supported by the Young core Teachers program of Henan province under
grant number 2017GGJS044.

\end{document}